\documentclass[11pt,a4paper]{article}

\usepackage[T1]{fontenc}
\usepackage{lmodern}
\usepackage{amsmath,amssymb,amsthm,mathtools}
\usepackage{booktabs}
\usepackage{array}
\usepackage{enumitem}
\usepackage{geometry}
\usepackage{indentfirst}
\usepackage{xcolor}
\usepackage{microtype}
\usepackage{hyperref}

\allowdisplaybreaks
\hypersetup{
  colorlinks=true,
  linkcolor=black,
  citecolor=black,
  urlcolor=black,
  pdftitle={Circular von Neumann--Jordan Constants in Complex Banach Spaces},
  pdfauthor={Haiying Li, Yijia Zhao and Yongjin Li}
}

\newtheorem{theorem}{Theorem}[section]
\newtheorem{proposition}[theorem]{Proposition}
\newtheorem{lemma}[theorem]{Lemma}
\newtheorem{corollary}[theorem]{Corollary}
\theoremstyle{definition}
\newtheorem{definition}[theorem]{Definition}
\newtheorem{example}[theorem]{Example}
\theoremstyle{remark}
\newtheorem{remark}[theorem]{Remark}

\let\C\relax
\newcommand{\C}{\mathbb C}
\newcommand{\R}{\mathbb R}
\newcommand{\T}{\mathbb T}
\newcommand{\CNJ}{C_{\mathrm{NJ}}}
\newcommand{\UC}{C_{\mathrm{NJ}}^{\circ}}
\newcommand{\LC}{c_{\mathrm{NJ}}^{\circ}}
\newcommand{\BC}{C_{\mathrm{NJ}}^{\mathbb C}}
\newcommand{\Um}{C_{\mathrm{NJ}}^{[m]}}
\newcommand{\Lm}{c_{\mathrm{NJ}}^{[m]}}
\newcommand{\dm}{\,dm}
\definecolor{RevisionBlue}{RGB}{0,70,190}
\newcommand{\rev}[1]{#1}
\newenvironment{revision}{\begingroup}{\endgroup}

\newcommand{\Synth}[1]{\rev{\mathcal{S}_{#1}}}
\newcommand{\SynthInv}[1]{\rev{\mathcal{S}_{#1}^{-1}}}

\title{\textbf{Circular von Neumann--Jordan Constants in Complex Banach Spaces}}

\author{Haiying Li, Yijia Zhao\\[2mm]
\small  School of Mathematics and Statistics,\\
\small  Henan Normal University,\\
\small  Xinxiang 453007\\
\small \texttt{lihaiying@htu.edu.cn(Haiying Li);zhaoyijia20259@stu.htu.edu.cn(Yijia Zhao)}\\[2mm]
Yongjin Li\\[2mm]
\small Department of Mathematics, Sun Yat-sen University,\\
\small Guangzhou 510275, Guangdong, China\\
\small \texttt{stslyj@mail.sysu.edu.cn}}
\date{}

\begin{document}
\maketitle

\begin{abstract}
Replacing \(y\) by \(ty\) does not refine the classical von
Neumann--Jordan constant; we instead average over unimodular phases
before optimization to define upper, lower, and balanced circular
constants for complex Banach spaces. The upper and lower constants
are the squared norm and squared minimum modulus of a Fourier
synthesis operator, linking them to optimal two-term Steinhaus type
and cotype constants. We compare them with Clarkson's constant and
characterize complex Hilbert spaces by the upper constant being one.
Root-of-unity averages provide finite-phase approximations with
explicit dimension-free errors. We compute exact one-sided values on \(L_p\)-spaces and Schatten
classes, the upper constant on \(\ell_1^n\), and all three constants
on \(\ell_\infty^n\), \(c_0\), \(\ell_\infty\), and \(C(K)\).
We also obtain quantitative comparisons with complex convexity and
smoothness moduli. Complex interpolation
improves the dimension-free upper bound for finite-dimensional
\(\ell_p\)-spaces, \(2<p<\infty\). A James constant estimate yields
the sharp uniform-nonsquareness threshold \(1+2/\pi\) and a sufficient
condition for the fixed point property of nonexpansive mappings.
\end{abstract}



\section{Introduction}

Clarkson's constant detects the failure of the real-symmetric
parallelogram identity.  Complex uniform convexity, Steinhaus
randomization, and PL-convexity are instead inherently circular; their
natural quadratic model averages a full orbit of unimodular phases.
This distinction motivates a complex invariant based on circular
\(L_2\)-means rather than on an additional scalar parameter.

For a real or complex Banach space, Clarkson's von Neumann--Jordan
constant is
\[
 \CNJ(X)=
 \sup_{(x,y)\ne(0,0)}
 \frac{\|x+y\|^2+\|x-y\|^2}
 {2(\|x\|^2+\|y\|^2)}.
 \tag{1.1}\label{eq:classical}
\]
It satisfies \(1\leq\CNJ(X)\leq2\), with equality
\(\CNJ(X)=1\) precisely for Hilbert spaces
\cite{Clarkson1937,JordanVN1935}.  Replacing \(y\) by \(ty\),
\(|t|\leq1\), does not refine \eqref{eq:classical}: the resulting
supremum is again \(\CNJ(X)\), as shown in Proposition
\ref{prop:collapse}.

Let \(\T=\{\zeta\in\C:|\zeta|=1\}\) and let \(m\) be normalized Haar
measure.  For a complex Banach space \(X\), set
\[
 M_X(x,y)^2=\int_{\T}\|x+\zeta y\|^2\dm(\zeta).
 \tag{1.2}\label{eq:mean}
\]
The upper and lower circular constants are the supremum and infimum of
\(M_X(x,y)^2/(\|x\|^2+\|y\|^2)\); the balanced constant is the maximum
of the upper constant and the reciprocal of the lower one.  Precise
definitions appear in Section~2.

Complex convexity and PL-convexity provide the analytic background for
\eqref{eq:mean} \cite{Globevnik1975,DavisGarlingTomczak1984,
Dilworth1986}.  Classical and modified von Neumann--Jordan constants,
including recent skew and multi-parameter variants, optimize over
finite configurations and are closely connected with the James
constant and normal structure \cite{KatoTakahashi1997,
AlonsoMartinPapini2008,CiesielskiPluciennik2020,FuLi2023,
NiLiuZhou2025,WanEtAl2026}.  Our
construction differs by averaging the phase before optimizing.  The
resulting Fourier operator permits duality and complex interpolation
arguments not visible in a free-parameter supremum.

The three principal results are the following.  They are proved in
Sections~4, 8, and 9, respectively.

\begin{revision}
\begin{theorem}\label{thm:intro-Hilbert}
For every nonzero complex Banach space \(X\),
\[
 \begin{gathered}
 1\leq\BC(X)\leq\CNJ(X)\leq2,\\
 \UC(X)=1\Longleftrightarrow\BC(X)=1
 \Longleftrightarrow X\text{ is a Hilbert space}.
 \end{gathered}
 \tag{1.3}\label{eq:intro-Hilbert}
\]
\end{theorem}
\end{revision}

\begin{theorem}
\label{thm:intro-fixed-point}
If \(\UC(X)<1+2/\pi\), then \(X\) is uniformly nonsquare and every
nonexpansive self-map of a nonempty bounded closed convex subset of
\(X\) has a fixed point.  More precisely, with
\(\gamma=1+2/\pi-\UC(X)>0\),
\[
 J(X)\leq2-\frac{\gamma^2}{32}.
 \tag{1.4}\label{eq:intro-James-gap}
\]
\rev{The threshold \(1+2/\pi\) for uniform nonsquareness is optimal.}
\end{theorem}

\begin{theorem}
\label{thm:intro-interpolation}
If \(n\geq2\) and \(2\leq p\leq\infty\), then
\[
 \UC(\ell_p^n)
 \leq
 \left(1+\frac n\pi\sin\frac\pi n\right)^{1-2/p}.
 \tag{1.5}\label{eq:intro-interpolation}
\]
For finite \(n\) and \(p>2\), this improves the classical
dimension-free estimate \(2^{1-2/p}\).
\end{theorem}

In addition, root-of-unity averages converge to the circular constants
with a dimension-free rate.  \rev{We compute all three constants on
Hilbert spaces, \(\ell_\infty^n\), \(c_0\), \(\ell_\infty\), and
\(C(K)\), the upper constant on \(\ell_1^n\), and one-sided values
on the \(L_p\) and Schatten scales.  The upper constant and the
reciprocal of the lower constant are the squared optimal two-term
Steinhaus type and cotype constants, respectively.}  Section~2 introduces
the circular constants, their operator representation, and their
relation to two-term Steinhaus type and cotype. Section~3 develops
finite-phase approximations. Section~4 establishes the Hilbert space
characterization and duality inequalities, while Section~5 relates
the constants to complex convexity and smoothness moduli.
Section~6 records stability and Banach--Mazur estimates.
Section~7 presents exact values and bounds in classical spaces. Section~8 provides the sharp uniform-nonsquareness
threshold \(1+2/\pi\) and a sufficient condition for the fixed point property of nonexpansive mappings.
Section~9 gives complex interpolation and dimension-sensitive bounds. Section~10 provides numerical checks of the preceding results.
Section~11 gives further problems.

\rev{All Banach spaces below are complex and nonzero.}  The unit sphere and closed unit
ball of \(X\) are denoted by \(S_X\) and \(B_X\).

\section{The circular constants and their operator meaning}

We first record the collapse of the disc-parameter definition in a
form that can be cited independently.

\begin{proposition}\label{prop:collapse}
For a complex Banach space \(X\), put
\[
 \widehat C_{\mathrm{NJ}}^{\mathbb C}(X)=
 \sup_{\substack{(x,y)\ne(0,0)\\ |t|\leq1}}
 \frac{\|x+ty\|^2+\|x-ty\|^2}
 {2(\|x\|^2+\|y\|^2)}.
\]
Then
\[
 \widehat C_{\mathrm{NJ}}^{\mathbb C}(X)=\CNJ(X).
\]
For a fixed \(t\), the corresponding supremum depends only on
\(|t|\), and every unimodular \(t\) gives \(\CNJ(X)\).
\end{proposition}

\begin{proof}
Fix \((x,y)\ne(0,0)\) and \(|t|\leq1\), and set \(z=ty\).
\rev{If \((x,z)=(0,0)\), the original quotient is zero and is at most
\(\CNJ(X)\).  Otherwise, its comparison with the quotient of
\((x,z)\) is legitimate.}  Since \(\|z\|\leq\|y\|\),
\[
 \frac{\|x+ty\|^2+\|x-ty\|^2}
 {2(\|x\|^2+\|y\|^2)}
 \leq
 \frac{\|x+z\|^2+\|x-z\|^2}
 {2(\|x\|^2+\|z\|^2)}
 \leq\CNJ(X).
\]
Taking \(t=1\) gives the reverse inequality.  If
\(t=re^{i\alpha}\), replacing \(y\) by \(e^{i\alpha}y\) shows that
the fixed-\(t\) supremum depends only on \(r=|t|\).  For \(r=1\) the
replacement is an isometry, so the fixed-\(t\) supremum is
\(\CNJ(X)\).
\end{proof}

\begin{definition}\label{def:circular}
For \(x,y\in X\), set
\[
 M_X(x,y)^2=\int_{\T}\|x+\zeta y\|^2\dm(\zeta).
 \tag{2.1}\label{eq:M-def}
\]
\rev{Here \(M_X(x,y)\) denotes the nonnegative square root of this
integral.}  The upper circular, lower circular and balanced complex
von Neumann--Jordan constants are respectively
\[
 \UC(X)=\sup_{(x,y)\ne(0,0)}
 \frac{M_X(x,y)^2}{\|x\|^2+\|y\|^2},
 \tag{2.2}\label{eq:upper-def}
\]
\[
 \LC(X)=\inf_{(x,y)\ne(0,0)}
 \frac{M_X(x,y)^2}{\|x\|^2+\|y\|^2},
 \tag{2.3}\label{eq:lower-def}
\]
and
\[
 \BC(X)=\max\left\{\UC(X),\frac1{\LC(X)}\right\}.
 \tag{2.4}\label{eq:balanced-def}
\]
\end{definition}

\rev{The next proposition gives the universal bounds; in particular, the lower constant is always positive.}

\begin{proposition}\label{prop:range}
For every complex Banach space \(X\),
\[
 \frac12\leq\LC(X)\leq1\leq\UC(X)\leq2.
 \tag{2.5}\label{eq:range}
\]
If \(H\) is a complex Hilbert space, then
\[
 \UC(H)=\LC(H)=\BC(H)=1.
 \tag{2.6}\label{eq:Hilbert-easy}
\]
\end{proposition}

\begin{proof}
The triangle inequality gives
\[
 M_X(x,y)^2
 \leq(\|x\|+\|y\|)^2
 \leq2(\|x\|^2+\|y\|^2).
\]
Taking \(y=0\) shows both \(\UC(X)\geq1\) and \(\LC(X)\leq1\).
Since
\[
 x=\int_{\T}(x+\zeta y)\dm(\zeta),\qquad
 y=\int_{\T}\overline\zeta(x+\zeta y)\dm(\zeta),
\]
Jensen's inequality yields
\[
 M_X(x,y)^2\geq\max\{\|x\|^2,\|y\|^2\}
 \geq\frac{\|x\|^2+\|y\|^2}{2}.
\]
This proves \eqref{eq:range}.

If \(H\) is a complex Hilbert space, then integration of the inner-product expansion
gives
\[
 M_H(x,y)^2=\|x\|^2+\|y\|^2,
\]
because \(\int_{\T}\zeta\dm(\zeta)=0\).
\end{proof}

\begin{example}\label{ex:endpoints}
\begin{enumerate}[label=\textup{(\roman*)}]
\item In complex \(\ell_1^2\), take the coordinate vectors
\(e_1,e_2\).  Since
\(\|e_1+\zeta e_2\|_1=2\) for every \(\zeta\in\T\), their quotient
in \eqref{eq:upper-def} equals \(2\).  Hence
\(\UC(\ell_1^2)=2\).
\item In complex \(\ell_\infty^2\), the same pair satisfies
\(\|e_1+\zeta e_2\|_\infty=1\).  Its lower quotient is \(1/2\), so
\(\LC(\ell_\infty^2)=1/2\).
\end{enumerate}
\end{example}

The constants have a direct operator interpretation.  Equip
\(X\oplus X\) with
\[
 \|(x,y)\|_2=(\|x\|^2+\|y\|^2)^{1/2},
\]
and define
\[
 \Synth{X}:X\oplus_2X\longrightarrow L_2(\T;X),\qquad
 (\Synth{X}(x,y))(\zeta)=x+\zeta y.
 \tag{2.7}\label{eq:S}
\]
For a bounded-below operator \(T\), write
\[
 j(T)=\inf_{\|u\|=1}\|Tu\|.
\]

\begin{proposition}\label{prop:operator}
The operator \(\Synth{X}\) is bounded below and
\[
 \UC(X)=\|\Synth{X}\|^2,\qquad
 \LC(X)=j(\Synth{X})^2.
 \tag{2.8}\label{eq:operator}
\]
Consequently,
\[
 \BC(X)=
 \max\{\|\Synth{X}\|^2,\|\SynthInv{X}:\Synth{X}(X\oplus_2X)\to
 X\oplus_2X\|^2\}.
 \tag{2.9}\label{eq:distortion}
\]
\end{proposition}

\begin{proof}
For every \((x,y)\in X\oplus_2X\), the function
\(\zeta\mapsto x+\zeta y\) is continuous on \(\T\), and the
Bochner \(L_2\)-norm gives
\[
 \|\Synth{X}(x,y)\|_{L_2(\T;X)}^2
 =\int_{\T}\|x+\zeta y\|^2\dm(\zeta)
 =M_X(x,y)^2.
\]
Since \(\|(x,y)\|_2^2=\|x\|^2+\|y\|^2\), homogeneity yields
\[
 \begin{aligned}
 \|\Synth{X}\|^2
 &=\sup_{\|x\|^2+\|y\|^2=1}M_X(x,y)^2
 =\sup_{(x,y)\ne(0,0)}
   \frac{M_X(x,y)^2}{\|x\|^2+\|y\|^2}
 =\UC(X),\\
 j(\Synth{X})^2
 &=\inf_{\|x\|^2+\|y\|^2=1}M_X(x,y)^2
 =\inf_{(x,y)\ne(0,0)}
   \frac{M_X(x,y)^2}{\|x\|^2+\|y\|^2}
 =\LC(X).
 \end{aligned}
\]
Proposition \ref{prop:range} therefore implies
\(\|\Synth{X}\|\leq\sqrt2\) and \(j(\Synth{X})\geq1/\sqrt2\).
Thus \(\Synth{X}\) is bounded and bounded below, hence injective with
closed range \(R_X=\Synth{X}(X\oplus_2X)\).  Its inverse on \(R_X\)
satisfies
\[
 \begin{aligned}
 \|\SynthInv{X}:R_X\to X\oplus_2X\|
 &=\sup_{0\ne u\in X\oplus_2X}
   \frac{\|u\|_2}{\|\Synth{X}u\|_{L_2(\T;X)}}\\
 &=\left(\inf_{0\ne u\in X\oplus_2X}
   \frac{\|\Synth{X}u\|_{L_2(\T;X)}}{\|u\|_2}\right)^{-1}
 =\frac1{j(\Synth{X})}.
 \end{aligned}
\]
Substituting these identities into
\(\BC(X)=\max\{\UC(X),\LC(X)^{-1}\}\) proves
\eqref{eq:distortion}.
\end{proof}

\begin{remark}
\rev{Formula \eqref{eq:distortion} expresses \(\BC(X)\) as the maximum
of the squared norm of the synthesis operator and the squared norm
of its inverse on its range.}  A Hilbert space makes the canonical
Fourier column \((1,\zeta)\) an isometry.
\end{remark}

The same operator has a probabilistic interpretation, in the spirit of
two-point type and cotype inequalities \cite{KatoTakahashiType1997}.
Let
\(\varepsilon_1,\varepsilon_2\) be independent signs and let
\(\zeta_1,\zeta_2\) be independent Steinhaus variables, that is,
independent random variables uniformly distributed on \(\T\).  Define
the optimal two-term constants by
\[
 T_{2,2}^{\T}(X)^2=
 \sup_{(x,y)\ne(0,0)}
 \frac{\mathbb E\|\zeta_1x+\zeta_2y\|^2}
 {\|x\|^2+\|y\|^2},
 \tag{2.10}\label{eq:Steinhaus-type}
\]
and
\[
 C_{2,2}^{\T}(X)^2=
 \sup_{(x,y)\ne(0,0)}
 \frac{\|x\|^2+\|y\|^2}
 {\mathbb E\|\zeta_1x+\zeta_2y\|^2}.
 \tag{2.11}\label{eq:Steinhaus-cotype}
\]

\begin{proposition}
\label{prop:Steinhaus}
For every complex Banach space \(X\),
\[
 T_{2,2}^{\T}(X)^2=\UC(X),\qquad
 C_{2,2}^{\T}(X)^2=\frac1{\LC(X)}.
 \tag{2.12}\label{eq:Steinhaus-identities}
\]
The corresponding two-term Rademacher type constant has square
\(\CNJ(X)\).  Hence Corollary \ref{cor:comparison} below compares
Steinhaus and Rademacher randomization at the exact two-vector level.
\end{proposition}

\begin{proof}
Fix \(x,y\in X\).  Since \(|\zeta_1|=1\), absolute homogeneity
of the complex norm gives
\[
 \|\zeta_1x+\zeta_2y\|
 =\|x+\overline{\zeta_1}\zeta_2y\|.
\]
Conditionally on \(\zeta_1\), independence and rotation invariance
of Haar measure imply that \(\overline{\zeta_1}\zeta_2\) is
uniformly distributed on \(\T\).  Consequently,
\[
 \mathbb E\|\zeta_1x+\zeta_2y\|^2
 =\int_{\T}\|x+\zeta y\|^2\dm(\zeta)
 =M_X(x,y)^2.
\]
Inserting this equality into \eqref{eq:Steinhaus-type} yields
\[
 T_{2,2}^{\T}(X)^2
 =\sup_{(x,y)\ne(0,0)}
   \frac{M_X(x,y)^2}{\|x\|^2+\|y\|^2}
 =\UC(X).
\]
For the cotype constant, Proposition \ref{prop:range} ensures that
the quotients defining \(\LC(X)\) are bounded below by \(1/2\).
Taking reciprocals therefore gives
\[
 \begin{aligned}
 C_{2,2}^{\T}(X)^2
 &=\sup_{(x,y)\ne(0,0)}
   \frac{\|x\|^2+\|y\|^2}{M_X(x,y)^2}\\
 &=\left(\inf_{(x,y)\ne(0,0)}
   \frac{M_X(x,y)^2}{\|x\|^2+\|y\|^2}\right)^{-1}
 =\frac1{\LC(X)}.
 \end{aligned}
\]
Finally, the four pairs of independent Rademacher signs are
equiprobable, so
\[
 \begin{aligned}
 \mathbb E\|\varepsilon_1x+\varepsilon_2y\|^2
 &=\frac14\bigl(\|x+y\|^2+\|x-y\|^2
              +\|-x+y\|^2+\|-x-y\|^2\bigr)\\
 &=\frac{\|x+y\|^2+\|x-y\|^2}{2}.
 \end{aligned}
\]
Dividing by \(\|x\|^2+\|y\|^2\) and taking the supremum over
nonzero pairs gives \(\CNJ(X)\), by \eqref{eq:classical}.
\end{proof}

\medskip
\noindent\textbf{Notation glossary and comparison.}
The following table fixes the notation used in the remainder of the
paper and separates the operations being measured.
\begin{center}
\small
\begin{tabular}{@{}c p{0.36\linewidth} p{0.38\linewidth}@{}}
\toprule
quantity & operation & geometric information\\
\midrule
\(\CNJ(X)\) & average over the two signs &
 upper parallelogram defect\\
\(\UC(X)\) & upper full-phase mean &
 circular expansion; Steinhaus type\\
\(\LC(X)\) & lower full-phase mean &
 circular contraction; reciprocal cotype\\
\(\BC(X)\) & maximum of the two distortions&
 two-sided Hilbertian defect\\
\(J(X)\) & \rev{supremum of the two-sign minimum} &
 uniform nonsquareness\\
\(\Um(X),\Lm(X)\) & mean over \(m\)-th roots of unity &
 finite-phase approximation\\
\bottomrule
\end{tabular}
\end{center}

\section{From finitely many phases to the circle}

Throughout this section, \(X\) is a nonzero complex Banach space,
\rev{and all phase counts are integers at least two.}  Let
\[
 \omega_m=e^{2\pi i/m}\qquad(m\geq2).
\]
Root-of-unity operator inequalities provide a related finite Fourier
framework; see \cite{BhatiaKittaneh2004}.
Define the finite-phase upper and lower constants by
\[
 \Um(X)=
 \sup_{(x,y)\ne(0,0)}
 \frac{\displaystyle\frac1m\sum_{k=0}^{m-1}
 \|x+\omega_m^k y\|^2}
 {\|x\|^2+\|y\|^2},
 \tag{3.1}\label{eq:Um}
\]
\[
 \Lm(X)=
 \inf_{(x,y)\ne(0,0)}
 \frac{\displaystyle\frac1m\sum_{k=0}^{m-1}
 \|x+\omega_m^k y\|^2}
 {\|x\|^2+\|y\|^2}.
 \tag{3.2}\label{eq:Lm}
\]
The case \(m=2\) is exactly the pair of phases \(1,-1\).

\begin{proposition}\label{prop:m2}
For every complex Banach space \(X\),
\[
 C_{\mathrm{NJ}}^{[2]}(X)=\CNJ(X),\qquad
 c_{\mathrm{NJ}}^{[2]}(X)=\frac1{\CNJ(X)}.
 \tag{3.3}\label{eq:m2}
\]
\end{proposition}

\begin{proof}
For \(m=2\), the upper quotient reduces precisely to Clarkson's
definition \eqref{eq:classical}.  Let
\[
 R(x,y)=
 \frac{\|x+y\|^2+\|x-y\|^2}
 {2(\|x\|^2+\|y\|^2)}.
\]
With \(u=x+y\) and \(v=x-y\), a direct calculation gives
\[
 R(u,v)=\frac1{R(x,y)}.
\]
The map \((x,y)\mapsto(u,v)\) is bijective, so
\(\inf R=1/\sup R\).
\end{proof}

\begin{theorem}\label{thm:refinement}
For integers \(q,m\geq2\),
\[
 C_{\mathrm{NJ}}^{[qm]}(X)\leq C_{\mathrm{NJ}}^{[m]}(X),
 \qquad
 c_{\mathrm{NJ}}^{[qm]}(X)\geq c_{\mathrm{NJ}}^{[m]}(X).
 \tag{3.4}\label{eq:divisibility}
\]
Moreover, for every \(m\geq2\),
\[
 \UC(X)\leq C_{\mathrm{NJ}}^{[m]}(X),\qquad
 \LC(X)\geq c_{\mathrm{NJ}}^{[m]}(X).
 \tag{3.5}\label{eq:circle-better}
\]
\end{theorem}

\begin{proof}
For \((x,y)\ne(0,0)\), write
\[
 R_m(x,y)=\frac{\displaystyle\frac1m\sum_{k=0}^{m-1}
 \|x+\omega_m^k y\|^2}{\|x\|^2+\|y\|^2}.
\]
For fixed \(x,y\), the function \(\alpha\mapsto R_m(x,\alpha y)\)
is continuous on \(\T\).  Moreover, the triangle inequality gives
\[
 0\leq R_m(x,\alpha y)
 \leq\frac{(\|x\|+\|y\|)^2}{\|x\|^2+\|y\|^2}\leq2.
\]
Thus all averages below are well defined.  For every \(\alpha\in\T\),
the pair \((x,\alpha y)\) is nonzero and has the same denominator, so
\[
 \Lm(X)\leq R_m(x,\alpha y)\leq\Um(X).
\]
Set \(\alpha_r=\omega_{qm}^{r}\) for \(0\leq r<q\).
Every index \(0\leq j<qm\) has a unique representation
\(j=r+qk\), with \(0\leq r<q\) and \(0\leq k<m\), and
\(\omega_{qm}^{r+qk}=\alpha_r\omega_m^k\).  Consequently,
\[
 R_{qm}(x,y)=\frac1q\sum_{r=0}^{q-1}R_m(x,\alpha_r y),
 \qquad
 \Lm(X)\leq R_{qm}(x,y)\leq\Um(X).
\]
Taking the infimum and supremum over all nonzero pairs proves
\eqref{eq:divisibility}.

For the circular comparison, rotation invariance of normalized Haar
measure gives, for each \(k\),
\[
 \int_{\T}\|x+\alpha\omega_m^k y\|^2\dm(\alpha)
 =\int_{\T}\|x+\alpha y\|^2\dm(\alpha)
 =M_X(x,y)^2.
\]
Interchanging the finite sum with the integral and dividing by
\(\|x\|^2+\|y\|^2\) yields
\[
 \frac{M_X(x,y)^2}{\|x\|^2+\|y\|^2}
 =\int_{\T}R_m(x,\alpha y)\dm(\alpha).
\]
Since Haar measure is normalized to have total mass one, this
integral lies between \(\Lm(X)\) and \(\Um(X)\).  Taking the infimum and supremum over nonzero pairs
now gives \eqref{eq:circle-better}.
\end{proof}

We next give a uniform error estimate.  It is useful both
conceptually and numerically.

\begin{lemma}\label{lem:Lipschitz}
If \(\|x\|^2+\|y\|^2=1\), then
\[
 g_{x,y}(\theta)=\|x+e^{i\theta}y\|^2
\]
is Lipschitz with constant at most \(1+\sqrt2\).
\end{lemma}

\begin{proof}
Put \(a=\|x\|\) and \(b=\|y\|\), so that
\(a^2+b^2=1\).  For \(\theta,\varphi\in\R\), set
\(A=\|x+e^{i\theta}y\|\) and
\(B=\|x+e^{i\varphi}y\|\).  The reverse triangle inequality and
\[
 |e^{i\theta}-e^{i\varphi}|
 =2\left|\sin\frac{\theta-\varphi}{2}\right|
 \leq|\theta-\varphi|
\]
imply \(|A-B|\leq b|\theta-\varphi|\).  The triangle
inequality also gives \(A+B\leq2(a+b)\).  Hence
\[
 |g_{x,y}(\theta)-g_{x,y}(\varphi)|
 =|A-B|(A+B)
 \leq2b(a+b)|\theta-\varphi|.
\]
To bound the coefficient, write \(a=\cos t\), \(b=\sin t\)
with \(t\in[0,\pi/2]\).  Then
\[
 2b(a+b)=1+\sin(2t)-\cos(2t)
 =1+\sqrt2\sin(2t-\pi/4)\leq1+\sqrt2.
\]
This proves the asserted Lipschitz bound, including the cases
\(a=0\) or \(b=0\).
\end{proof}

\begin{theorem}\label{thm:limit}
For every integer \(m\geq2\),
\[
 \left|C_{\mathrm{NJ}}^{[m]}(X)-\UC(X)\right|
 \leq\frac{(1+\sqrt2)\pi}{2m},
 \tag{3.6}\label{eq:upper-rate}
\]
and
\[
 \left|c_{\mathrm{NJ}}^{[m]}(X)-\LC(X)\right|
 \leq\frac{(1+\sqrt2)\pi}{2m}.
 \tag{3.7}\label{eq:lower-rate}
\]
In particular, both finite-phase constants converge to their
circular counterparts.
\end{theorem}

\begin{proof}
For a nonzero pair \((x,y)\), set
\(r=(\|x\|^2+\|y\|^2)^{1/2}>0\).
Replacing \((x,y)\) by \((x/r,y/r)\) leaves both the finite-phase
and circular quotients unchanged.  Their suprema and infima can
therefore be taken over the nonempty set
\[
 \mathcal P=\{(x,y)\in X\times X:\|x\|^2+\|y\|^2=1\}.
\]
Fix \((x,y)\in\mathcal P\).  The function \(g=g_{x,y}\) is
continuous and \(2\pi\)-periodic, with
\(0\leq g(\theta)\leq(\|x\|+\|y\|)^2\leq2\).  Set
\[
 A_m(x,y)=\frac1m\sum_{k=0}^{m-1}g(\theta_k),\qquad
 B(x,y)=\frac1{2\pi}\int_0^{2\pi}g(\theta)\,d\theta,
 \qquad \theta_k=\frac{2\pi k}{m}.
\]
Let \(I_k=[\theta_k-\pi/m,\theta_k+\pi/m)\).  These intervals
partition \([-\pi/m,2\pi-\pi/m)\), a full period of \(g\).
Since each \(I_k\) has length \(2\pi/m\), periodicity gives
\[
 A_m(x,y)-B(x,y)
 =\frac1{2\pi}\sum_{k=0}^{m-1}
   \int_{I_k}\bigl(g(\theta_k)-g(\theta)\bigr)\,d\theta.
\]
Applying Lemma \ref{lem:Lipschitz} on each interval, we obtain
\[
 \begin{aligned}
 |A_m(x,y)-B(x,y)|
 &\leq\frac{1+\sqrt2}{2\pi}
       \sum_{k=0}^{m-1}\int_{I_k}|\theta-\theta_k|\,d\theta\\
 &=\frac{(1+\sqrt2)m}{2\pi}
       \int_{-\pi/m}^{\pi/m}|s|\,ds
 =\frac{(1+\sqrt2)\pi}{2m}=: \varepsilon_m.
 \end{aligned}
\]
The bound is uniform over \(\mathcal P\) and independent of
\(X\) and its dimension.  Thus
\(B-\varepsilon_m\leq A_m\leq B+\varepsilon_m\) throughout
\(\mathcal P\), and taking suprema and infima gives
\[
 \begin{aligned}
 \sup_{\mathcal P}B-\varepsilon_m
 &\leq\sup_{\mathcal P}A_m
 \leq\sup_{\mathcal P}B+\varepsilon_m,\\
 \inf_{\mathcal P}B-\varepsilon_m
 &\leq\inf_{\mathcal P}A_m
 \leq\inf_{\mathcal P}B+\varepsilon_m.
 \end{aligned}
\]
The respective suprema are \(\UC(X)\) and \(\Um(X)\), and
the respective infima are \(\LC(X)\) and \(\Lm(X)\).
These extrema need not be attained: the inequalities use only the
order-preserving properties of supremum and infimum.
This proves \eqref{eq:upper-rate} and \eqref{eq:lower-rate}.
Both convergence assertions follow from \(\varepsilon_m\to0\).
\end{proof}

\begin{corollary}
\label{cor:comparison}
For every complex Banach space \(X\),
\[
 \UC(X)\leq\CNJ(X),\qquad
 \LC(X)\geq\frac1{\CNJ(X)},
 \tag{3.8}\label{eq:comparison}
\]
and therefore
\[
 1\leq\BC(X)\leq\CNJ(X)\leq2.
 \tag{3.9}\label{eq:balanced-comparison}
\]
\end{corollary}

\begin{proof}
Use \eqref{eq:circle-better} with \(m=2\) and then Proposition
\ref{prop:m2}.
\end{proof}

\begin{remark}
The circular constant is not obtained by making the phase set larger
inside a supremum.  The order of operations is essential: the phases
are averaged first, and the vectors are optimized afterwards.  This
is why \(\UC(X)\leq\CNJ(X)\), whereas the naive intuition might
suggest the opposite inequality.
\end{remark}

\section{Hilbert spaces, duality and two-sided control}

The Hilbertian mean-square identity
\[
 M_X(x,y)^2=\|x\|^2+\|y\|^2
 \tag{4.1}\label{eq:mean-identity}
\]
\rev{characterizes Hilbert spaces.  In fact, its upper inequality with
constant one already forces the full identity, as proved in
Theorem~\ref{thm:Hilbert}.}

\begin{lemma}
\label{lem:Gaussian}
If \eqref{eq:mean-identity} holds for all \(x,y\in X\), then \(X\)
satisfies the parallelogram identity.
\end{lemma}

\begin{proof}
Fix \(n\geq1\) and \(x_1,\ldots,x_n\in X\), and let
\(\zeta_1,\ldots,\zeta_n\) be independent Steinhaus variables.
Put \(Z_0=0\), \(Z_k=\sum_{j=1}^k\zeta_jx_j\), and
\(\mathcal F_k=\sigma(\zeta_1,\ldots,\zeta_k)\), with
\(\mathcal F_0\) the trivial sigma-algebra.  The bound
\(\|Z_k\|\leq\sum_{j=1}^k\|x_j\|\) ensures integrability of
\(\|Z_k\|^2\).  Independence and \eqref{eq:mean-identity} give,
almost surely,
\[
 \begin{aligned}
 \mathbb E\bigl(\|Z_k\|^2\mid\mathcal F_{k-1}\bigr)
 &=\int_{\T}\|Z_{k-1}+\zeta x_k\|^2\dm(\zeta)\\
 &=\|Z_{k-1}\|^2+\|x_k\|^2.
 \end{aligned}
\]
Taking expectations and iterating from \(Z_0=0\) proves
\[
 \mathbb E\left\|\sum_{j=1}^n\zeta_jx_j\right\|^2
 =\sum_{j=1}^n\|x_j\|^2.
 \tag{4.2}\label{eq:Steinhaus}
\]

Now let \(g_j=(\xi_j+i\eta_j)/\sqrt2\), where all
\(\xi_j,\eta_j\) are independent real standard Gaussian variables.
Each \(g_j\) has planar density \(\pi^{-1}e^{-|z|^2}\).
Writing \(g_j=r_j\omega_j\), where \(r_j=|g_j|\), polar
coordinates give the product law
\(2r e^{-r^2}\,dr\dm(\omega)\) for \(r>0\) and \(\omega\in\T\).
Thus the \(\omega_j\) are independent Steinhaus variables,
independent of all the radii, and \(\mathbb E r_j^2=1\).
The phases can be assigned arbitrarily on the null events
\(g_j=0\).  Moreover,
\[
 \left\|\sum_{j=1}^n g_jx_j\right\|^2
 \leq\left(\sum_{j=1}^n|g_j|^2\right)
       \left(\sum_{j=1}^n\|x_j\|^2\right),
\]
so the squared norm is integrable.  Applying
\eqref{eq:Steinhaus} conditionally on the radii yields
\[
 \mathbb E\left(\left\|\sum_{j=1}^n g_jx_j\right\|^2
       \mathrel{\Big|}r_1,\ldots,r_n\right)
 =\sum_{j=1}^n r_j^2\|x_j\|^2.
\]
The tower property therefore gives
\[
 \mathbb E\left\|\sum_{j=1}^n g_jx_j\right\|^2
 =\sum_{j=1}^n\|x_j\|^2.
 \tag{4.3}\label{eq:Gaussian}
\]

Finally, take \(n=2\) and put
\[
 h_1=\frac{g_1+g_2}{\sqrt2},\qquad
 h_2=\frac{g_1-g_2}{\sqrt2}.
\]
The defining matrix is unitary, so it preserves both
\(|z_1|^2+|z_2|^2\) and Lebesgue measure on \(\C^2\).
It therefore preserves the joint density
\(\pi^{-2}e^{-(|z_1|^2+|z_2|^2)}\), and
\((h_1,h_2)\) has the same distribution as \((g_1,g_2)\).
Using \eqref{eq:Gaussian} twice gives the parallelogram identity:
\begin{align*}
 \|x\|^2+\|y\|^2
 &=\mathbb E\|h_1x+h_2y\|^2\\
 &=\mathbb E\left\|g_1\frac{x+y}{\sqrt2}
                    +g_2\frac{x-y}{\sqrt2}\right\|^2\\
 &=\frac{\|x+y\|^2+\|x-y\|^2}{2}
 \qquad(x,y\in X).\qedhere
\end{align*}
\end{proof}

\begin{revision}
\begin{theorem}\label{thm:Hilbert}
For a nonzero complex Banach space \(X\), the following statements
are equivalent:
\begin{enumerate}[label=\textup{(\roman*)}]
\item the original norm of \(X\) is induced by a complex inner
product, so \(X\) is a Hilbert space;
\item \(\BC(X)=1\);
\item \(\UC(X)=\LC(X)=1\);
\item the mean-square identity \eqref{eq:mean-identity} holds for all
\(x,y\in X\);
\item \(\UC(X)=1\).
\end{enumerate}
\end{theorem}

\begin{proof}
Proposition~\ref{prop:range} gives (i)\(\Rightarrow\)(iii).
The definition of \(\BC\), together with
\(\LC(X)\leq1\leq\UC(X)\), gives
(ii)\(\Leftrightarrow\)(iii), while the definitions of the upper
and lower constants give (iii)\(\Leftrightarrow\)(iv).
Lemma~\ref{lem:Gaussian} and the Jordan--von Neumann theorem imply
(iv)\(\Rightarrow\)(i) \cite{JordanVN1935}; also,
(iii)\(\Rightarrow\)(v).  It remains to prove
(v)\(\Rightarrow\)(i).  We use a complex version of the contact-point
argument for Gaussian type two with constant one; compare
\cite[Proposition~3.13]{Yaroslavtsev2020}.  The details below derive
the required complex Gaussian estimate directly from (v).

\smallskip
\noindent\hspace*{2em}\textbf{\emph{Step 1: A two-vector Gaussian estimate.}}
Assume \(\UC(X)=1\).  Its supremum definition gives
\[
 \int_{\T}\|u+\zeta v\|^2\dm(\zeta)
 \leq\|u\|^2+\|v\|^2\qquad(u,v\in X),
\]
including \((u,v)=(0,0)\).  Let \(g_1,g_2\) be independent
standard complex Gaussian variables, each with density
\(\pi^{-1}e^{-|z|^2}\) on \(\C\).  Write
\(g_j=r_j\zeta_j\), where the radii have density
\(2r e^{-r^2}\) on \([0,\infty)\), the phases are independent
Steinhaus variables, and the radii and phases are independent.
In particular, \(\mathbb E r_j^2=1\) and \(\mathbb E g_j=0\).
For fixed radii, homogeneity and Haar invariance yield
\[
 \mathbb E_{\zeta_1,\zeta_2}
   \|r_1\zeta_1u+r_2\zeta_2v\|^2
 =M_X(r_1u,r_2v)^2
 \leq r_1^2\|u\|^2+r_2^2\|v\|^2.
\]
Integrating the radii gives
\[
 \mathbb E\|g_1u+g_2v\|^2\leq\|u\|^2+\|v\|^2.
 \tag{4.3a}\label{eq:Gaussian-upper-one}
\]
Tonelli's theorem applies to the nonnegative integrands, and the
bound
\(\|g_1u+g_2v\|^2\leq2|g_1|^2\|u\|^2+2|g_2|^2\|v\|^2\)
also verifies finiteness of every expectation used here.

\smallskip
\noindent\hspace*{2em}\textbf{\emph{Step 2: An enclosing Hermitian ellipsoid with spanning
contact points.}}
Fix a complex two-dimensional subspace \(E\subset X\), with its
inherited norm, and identify it with \(\C^2\).
Let \(B=\{z\in E:\|z\|\leq1\}\), and let \(B_2\) be the unit
ball for a coordinate Euclidean norm \(|\cdot|_0\).
Choose \(r,R>0\) with \(rB_2\subset B\subset RB_2\).
Consider the feasible Hermitian matrices
\[
 \mathcal H=\{H=H^*\geq0:z^*Hz\leq1\ \text{for all }z\in B\}.
\]
This set is closed, and \(rB_2\subset B\) implies
\(0\leq H\leq r^{-2}I\), so it is compact.
Moreover, \(R^{-2}I\) is feasible and has positive determinant.
Hence a maximizer \(H_0\) of \(\det H\) on \(\mathcal H\) is
positive definite.  The Hermitian norm
\[
 |z|=(z^*H_0z)^{1/2}
\]
satisfies \(|z|\leq\|z\|\) for all \(z\in E\).

We show that its contact set
\(K=\{z:|z|=\|z\|=1\}\) spans \(E\) over \(\C\).
The complex-linear change of coordinates \(w=H_0^{1/2}z\)
transforms \(B\) to \(B'\subset B_2\), the unit ball of
\(\|w\|'=\|H_0^{-1/2}w\|\).
A feasible matrix \(G\) in these coordinates corresponds to
\(H=H_0^{1/2}GH_0^{1/2}\), with
\(\det H=(\det H_0)\det G\).
Thus \(I\) is a determinant maximizer for \(B'\).
The new contact set is
\(K'=\{w\in B':|w|_0=1\}\).
It is nonempty: otherwise compactness would give
\(\max_{B'}|w|_0^2<1\), so a scalar multiple of \(I\) larger than
\(I\) would be feasible.  Each \(w\in K'\) also satisfies \(\|w\|'=1\), since
otherwise it could be scaled up within
\(B'\), contradicting \(B'\subset B_2\).

Suppose \(V=\operatorname{span}_{\C}K'\) is a proper subspace.
Let \(W=V^\perp\), \(q=\dim_{\C}W\geq1\),
\(\delta=q/4\), and \(A=P_W-\delta I\), where \(P_W\) is the
Euclidean orthogonal projection.  Then
\[
 \operatorname{tr}A=q-2\delta=q/2>0,
 \qquad w^*Aw=-\delta\quad(w\in K').
\]
The relatively open set
\(U=\{w\in B':w^*Aw<-\delta/2\}\) contains \(K'\).
If \(B'\setminus U\ne\varnothing\), compactness gives
\[
 \eta=1-\max_{w\in B'\setminus U}|w|_0^2>0.
\]
Set \(L=\max_{w\in B'}|w^*Aw|<\infty\).
For all sufficiently small \(\varepsilon>0\),
\(I+\varepsilon A\) is positive definite.
On \(U\) its quadratic form is at most
\(|w|_0^2-\varepsilon\delta/2\leq1\); on
\(B'\setminus U\) it is at most
\(1-\eta+\varepsilon L\leq1\) if \(\varepsilon L\leq\eta\).
If the complement is empty, only the first estimate is needed.
Thus one sufficiently small positive \(\varepsilon\) makes
\(I+\varepsilon A\) feasible on all of \(B'\).  However,
\[
 \det(I+\varepsilon A)
 =(1+\varepsilon(1-\delta))^q(1-\varepsilon\delta)^{2-q},
 \qquad
 \left.\frac{d}{d\varepsilon}\det(I+\varepsilon A)\right|_{\varepsilon=0}
 =q/2>0.
\]
Its determinant is therefore greater than one for sufficiently
small positive \(\varepsilon\), a contradiction.
It follows that \(K'\), and hence \(K\), spans \(E\).

\smallskip
\noindent\hspace*{2em}\textbf{\emph{Step 3: Equality of the original and Hermitian norms.}}
Choose complex-linearly independent \(z_1,z_2\in K\), so
\(|z_j|=\|z_j\|=1\), and put \(Z=g_1z_1+g_2z_2\).
Independence and \(\mathbb E(g_1\overline{g_2})=0\) give
\(\mathbb E|Z|^2=|z_1|^2+|z_2|^2=2\), without requiring
\(z_1,z_2\) to be orthogonal.  By \(|w|\leq\|w\|\) and
\eqref{eq:Gaussian-upper-one},
\[
 2=\mathbb E|Z|^2\leq\mathbb E\|Z\|^2
 \leq\|z_1\|^2+\|z_2\|^2=2.
\]
Consequently, the continuous nonnegative function
\(F(w)=\|w\|^2-|w|^2\) satisfies \(\mathbb EF(Z)=0\).
Since \(z_1,z_2\) form a complex basis of \(E\), the invertible
map \((a,b)\mapsto az_1+bz_2\) sends the everywhere positive
Gaussian density on \(\C^2\) to a distribution with full support
on \(E\).  If \(F(w_0)>0\), continuity gives a nonempty open
neighborhood \(O\) of \(w_0\) and a constant \(c>0\) such that
\(F\geq c\) on \(O\).  Then
\(\mathbb EF(Z)\geq c\mathbb P(Z\in O)>0\), a contradiction.
Thus \(F\equiv0\): the original norm on every complex
two-dimensional subspace is a Hermitian norm.

For arbitrary \(x,y\in X\), apply this conclusion to
\(E=\operatorname{span}_{\C}\{x,y\}\) when \(\dim_{\C}E=2\).
If its dimension is at most one, complex homogeneity gives the
parallelogram identity directly.  Hence in all cases
\[
 \|x+y\|^2+\|x-y\|^2=2\|x\|^2+2\|y\|^2.
\]
The Jordan--von Neumann theorem yields a complex inner product
inducing the original norm.  With the first variable linear, it is
\[
 \langle x,y\rangle
 =\frac14\bigl(\|x+y\|^2-\|x-y\|^2
   +i\|x+iy\|^2-i\|x-iy\|^2\bigr).
\]
Completeness of the original norm makes \(X\) a Hilbert space.
This proves (v)\(\Rightarrow\)(i) and completes the equivalences.
\end{proof}
\end{revision}

\rev{The next inequalities give quantitative dual control between
the lower constant of \(X\) and the upper constant of \(X^*\).}

\begin{theorem}\label{thm:dual}
For every nonzero complex Banach space \(X\),
\[
 \LC(X)\UC(X^*)\geq1,\qquad
 \LC(X^*)\UC(X)\geq1.
 \tag{4.4}\label{eq:dual}
\]
Consequently, if \(\UC(X)=\UC(X^*)=1\), then \(X\) is a Hilbert
space.
\end{theorem}

\begin{proof}
We use the complex-linear duality pairing \(f(x)\), with
\(f\in X^*\) and \(x\in X\).  All circle integrands below are
continuous, so the integrals are well defined.  Fix
\((x,y)\ne(0,0)\) and put \(D=\|x\|^2+\|y\|^2>0\).
By the Hahn--Banach theorem, choose \(x^*,y^*\in X^*\) satisfying
\[
 \begin{aligned}
 \|x^*\|=\frac{\|x\|}{\sqrt D},\qquad
 &x^*(x)=\frac{\|x\|^2}{\sqrt D},\\
 \|y^*\|=\frac{\|y\|}{\sqrt D},\qquad
 &y^*(y)=\frac{\|y\|^2}{\sqrt D},
 \end{aligned}
\]
where the corresponding functional is zero if its vector is zero.
Then \(\|x^*\|^2+\|y^*\|^2=1\).  Expanding the pairing and using
\(\int_{\T}\zeta\dm(\zeta)
=\int_{\T}\overline\zeta\dm(\zeta)=0\), we obtain
\[
 \begin{aligned}
 &\int_{\T}(x^*+\overline\zeta y^*)(x+\zeta y)\dm(\zeta)\\
 &\quad=\int_{\T}\bigl(x^*(x)+\zeta x^*(y)
              +\overline\zeta y^*(x)+y^*(y)\bigr)\dm(\zeta)
 =\sqrt D.
 \end{aligned}
\]
The triangle inequality for the scalar integral, the duality bound,
and Cauchy--Schwarz give
\[
 \begin{aligned}
 \sqrt D
 &=\left|\int_{\T}(x^*+\overline\zeta y^*)(x+\zeta y)\dm(\zeta)\right|\\
 &\leq\int_{\T}\|x^*+\overline\zeta y^*\|\,\|x+\zeta y\|\dm(\zeta)\\
 &\leq M_{X^*}(x^*,y^*)\,M_X(x,y)
 \leq\sqrt{\UC(X^*)}\,M_X(x,y).
 \end{aligned}
\]
Here Haar measure is invariant under
\(\zeta\mapsto\overline\zeta\), and the last inequality uses
the normalization of \((x^*,y^*)\).  The upper constant is
finite and at least one by Proposition \ref{prop:range}.
Squaring and rearranging gives
\[
 \frac{M_X(x,y)^2}{\|x\|^2+\|y\|^2}
 \geq\frac{1}{\UC(X^*)}.
\]
Taking the infimum over nonzero pairs proves
\(\LC(X)\UC(X^*)\geq1\).

For the second inequality, fix \((f,g)\ne(0,0)\) in
\(X^*\times X^*\), set \(D=\|f\|^2+\|g\|^2>0\), and take
\(0<\varepsilon<1\).  By the definition of the dual norm and a
choice of unimodular factors, there exist \(u,v\in S_X\) for which
\(f(u),g(v)\) are nonnegative real numbers and
\[
 f(u)\geq(1-\varepsilon)\|f\|,\qquad
 g(v)\geq(1-\varepsilon)\|g\|.
\]
If either functional is zero, its unit vector can be chosen
arbitrarily; the sphere is nonempty because \(X\ne\{0\}\).
Set
\[
 u_0=\frac{\|f\|}{\sqrt D}u,\qquad
 v_0=\frac{\|g\|}{\sqrt D}v,
 \qquad \|u_0\|^2+\|v_0\|^2=1.
\]
Integrating the expanded pairing cancels the mixed terms and gives
\[
 \int_{\T}(f+\zeta g)(u_0+\overline\zeta v_0)\dm(\zeta)
 =f(u_0)+g(v_0)\geq(1-\varepsilon)\sqrt D.
\]
Applying the same integral estimates and conjugation invariance,
we obtain
\[
 \begin{aligned}
 (1-\varepsilon)\sqrt D
 &\leq\left|\int_{\T}(f+\zeta g)(u_0+\overline\zeta v_0)\dm(\zeta)\right|\\
 &\leq M_{X^*}(f,g)\,M_X(u_0,v_0)
 \leq\sqrt{\UC(X)}\,M_{X^*}(f,g).
 \end{aligned}
\]
Squaring and rearranging yields
\[
 \frac{M_{X^*}(f,g)^2}{\|f\|^2+\|g\|^2}
 \geq\frac{(1-\varepsilon)^2}{\UC(X)}.
\]
Letting \(\varepsilon\downarrow0\) for this fixed pair and then
taking the infimum over all nonzero \((f,g)\) proves
\(\LC(X^*)\UC(X)\geq1\).  Thus the argument requires only
approximate norming vectors, rather than attainment of the dual norm.

If \(\UC(X)=\UC(X^*)=1\), the first inequality in
\eqref{eq:dual} gives \(\LC(X)\geq1\), whereas Proposition
\ref{prop:range} gives \(\LC(X)\leq1\).  Hence
\(\UC(X)=\LC(X)=1\), and Theorem \ref{thm:Hilbert} implies
that \(X\) is a Hilbert space.
\end{proof}

\begin{remark}
\rev{The duality inequalities provide general quantitative control.
At the extremal value \(\UC(X)=1\), however,
Theorem~\ref{thm:Hilbert} already implies the Hilbert structure and
therefore \(\LC(X)=\UC(X^*)=1\).  Thus the additional dual condition
in the preceding sufficient criterion is unnecessary at this endpoint.}
\end{remark}

\section{The connection with complex uniform convexity}

\rev{We relate the circular constants to moduli arising in complex
convexity and PL-convexity.}  Here PL refers to
\emph{plurisubharmonicity}: the norm
restricted to each complex affine line is subharmonic, so that
\[
 \|x\|\leq\int_{\T}\|x+\zeta y\|\dm(\zeta)
 \qquad(x,y\in X).
\]
This submean property underlies the terminology of PL-convexity
\cite{DavisGarlingTomczak1984,Dilworth1986,Pavlovic1991}.  In the two-vector
quadratic setting used here, the PL-convexity modulus
\(\delta_{2,X}^{\circ}(t)\) measures the smallest excess above
\(\|x\|=1\) of the circular \(L_2\)-mean of \(x+t\zeta y\),
over \(x,y\in S_X\); its PL-smoothness counterpart
\(\rho_{2,X}^{\circ}(t)\) measures the largest such excess.
Their precise definitions are given in \eqref{eq:rho} and
\eqref{eq:delta}.  It is important to distinguish these circular
averages from a maximum over a complex disc.

Throughout this section, \(X\) is a nonzero complex Banach space.
For \(t>0\), define the standard disc modulus
\[
 \Delta_X^{\mathbb C}(t)=
 \inf_{x,y\in S_X}
 \left[
 \sup_{|\lambda|\leq1}\|x+t\lambda y\|-1
 \right].
 \tag{5.1}\label{eq:disc-modulus}
\]
The space \(X\) is complex uniformly convex if and only if
\(\Delta_X^{\mathbb C}(t)>0\) for every \(t>0\); this is an
equivalent quantitative form of the definition in
\cite{Globevnik1975,DavisGarlingTomczak1984}.  By the maximum
principle for Banach-valued analytic functions
\cite{ThorpWhitley1967}, the supremum in
\eqref{eq:disc-modulus} may be taken over \(|\lambda|=1\).

For \(t\geq0\), the quadratic circular moduli are
\[
 \rho_{2,X}^{\circ}(t)=
 \sup_{x,y\in S_X}
 \left[
 \left(\int_{\T}\|x+t\zeta y\|^2\dm(\zeta)\right)^{1/2}-1
 \right],
 \tag{5.2}\label{eq:rho}
\]
\[
 \delta_{2,X}^{\circ}(t)=
 \inf_{x,y\in S_X}
 \left[
 \left(\int_{\T}\|x+t\zeta y\|^2\dm(\zeta)\right)^{1/2}-1
 \right].
 \tag{5.3}\label{eq:delta}
\]
The first is a two-vector quadratic PL-smoothness modulus; the second
is its PL-convexity companion.  \rev{The submean inequality and
Cauchy--Schwarz show that both moduli are nonnegative.}  In particular,
\(\rho_{2,X}^{\circ}(0)=\delta_{2,X}^{\circ}(0)=0\).

\begin{theorem}\label{thm:moduli}
For every nonzero complex Banach space \(X\),
\[
 \UC(X)=
 \sup_{t\geq0}
 \frac{(1+\rho_{2,X}^{\circ}(t))^2}{1+t^2},
 \tag{5.4}\label{eq:U-rho}
\]
\[
 \LC(X)=
 \inf_{t\geq0}
 \frac{(1+\delta_{2,X}^{\circ}(t))^2}{1+t^2}.
 \tag{5.5}\label{eq:L-delta}
\]
Moreover,
\[
 \Delta_X^{\mathbb C}(t)\geq
 \delta_{2,X}^{\circ}(t)\qquad(t>0).
 \tag{5.6}\label{eq:modulus-comparison}
\]
\end{theorem}

\begin{proof}
Consider the quotient in \eqref{eq:upper-def}--\eqref{eq:lower-def}
for \((x,y)\ne(0,0)\).  If exactly one coordinate is zero, then
\[
 \frac{M_X(x,y)^2}{\|x\|^2+\|y\|^2}=1.
\]
For \(x,y\ne0\), write
\[
 x=a u,\qquad y=b v,\qquad
 a=\|x\|>0,\quad b=\|y\|>0,\quad u,v\in S_X,
\]
and set \(t=b/a>0\).  Homogeneity gives
\[
 \begin{aligned}
 \frac{M_X(x,y)^2}{\|x\|^2+\|y\|^2}
 &=\frac{a^2\displaystyle\int_{\T}\|u+t\zeta v\|^2\dm(\zeta)}
         {a^2(1+t^2)}\\
 &=\frac{M_X(u,tv)^2}{1+t^2}.
 \end{aligned}
\]
Conversely, given \(t>0\) and \(u,v\in S_X\), the choice
\((x,y)=(u,tv)\) realizes the quotient on the right.
For \(t=0\), it equals \(1\), which also represents every
zero-coordinate case.  Thus the quotients obtained from all nonzero
pairs are exactly those parametrized by \(t\geq0\) and
\(u,v\in S_X\).

For each fixed \(t\geq0\), the set of values \(M_X(u,tv)\)
is nonempty and contained in \([0,1+t]\).  Squaring is continuous
and increasing on this interval, so the definitions of the moduli give
\[
 \begin{aligned}
 \sup_{u,v\in S_X}\frac{M_X(u,tv)^2}{1+t^2}
   &=\frac{\bigl(1+\rho_{2,X}^{\circ}(t)\bigr)^2}{1+t^2},\\
 \inf_{u,v\in S_X}\frac{M_X(u,tv)^2}{1+t^2}
   &=\frac{\bigl(1+\delta_{2,X}^{\circ}(t)\bigr)^2}{1+t^2}.
 \end{aligned}
\]
Taking the supremum over \(t\geq0\) in the first identity and
the infimum over \(t\geq0\) in the second, the preceding
parametrization proves \eqref{eq:U-rho} and
\eqref{eq:L-delta}.  These operations do not require attainment
of any extremum.

Finally, for every \(t>0\) and \(u,v\in S_X\), the probability
normalization \(m(\T)=1\) implies
\[
 M_X(u,tv)
 \leq\sup_{\zeta\in\T}\|u+t\zeta v\|
 \leq\sup_{|\lambda|\leq1}\|u+t\lambda v\|.
\]
Subtracting \(1\) and taking the infimum over the same unit
vector pairs yields \eqref{eq:modulus-comparison}.
\end{proof}

\begin{corollary}\label{cor:moduli}
For every nonzero complex Banach space \(X\) and every \(t\geq0\),
\[
 \rho_{2,X}^{\circ}(t)
 \leq\sqrt{\UC(X)(1+t^2)}-1,
 \tag{5.7}\label{eq:rho-bound}
\]
\[
 \delta_{2,X}^{\circ}(t)
 \geq\sqrt{\LC(X)(1+t^2)}-1.
 \tag{5.8}\label{eq:delta-bound}
\]
In particular, for every \(t>0\), the right-hand side of
\eqref{eq:delta-bound} is also a lower bound for
\(\Delta_X^{\mathbb C}(t)\).
\end{corollary}

\begin{proof}
Fix \(t\geq0\) and \(u,v\in S_X\).  The pair \((u,tv)\)
is nonzero even when \(t=0\), and
\(\|u\|^2+\|tv\|^2=1+t^2\).  Hence the definitions of the
upper and lower circular constants give
\[
 \LC(X)(1+t^2)\leq M_X(u,tv)^2\leq\UC(X)(1+t^2).
\]
Taking nonnegative square roots and subtracting \(1\), we obtain
\[
 \sqrt{\LC(X)(1+t^2)}-1
 \leq M_X(u,tv)-1
 \leq\sqrt{\UC(X)(1+t^2)}-1.
\]
Both bounds are independent of \(u,v\).  Taking the supremum
of the middle term in the upper estimate gives
\eqref{eq:rho-bound}; taking its infimum in the lower estimate
gives \eqref{eq:delta-bound}.  At \(t=0\), these inequalities
reduce to \(0\leq\sqrt{\UC(X)}-1\) and
\(0\geq\sqrt{\LC(X)}-1\), in agreement with
Proposition~\ref{prop:range}.

For every \(t>0\), \eqref{eq:modulus-comparison} further yields
\[
 \Delta_X^{\mathbb C}(t)
 \geq\delta_{2,X}^{\circ}(t)
 \geq\sqrt{\LC(X)(1+t^2)}-1.
\]
This proves the final assertion.
\end{proof}

\begin{example}\label{ex:Hilbert-moduli}
For a complex Hilbert space \(H\), \rev{the mean-square identity
\eqref{eq:mean-identity} gives}
\[
 \rho_{2,H}^{\circ}(t)=\delta_{2,H}^{\circ}(t)
 =\sqrt{1+t^2}-1.
\]
Thus \(\BC(X)\) measures the global two-sided departure of the
quadratic circular moduli from their Hilbertian benchmark.
\end{example}

\begin{remark}
The use of \(\sup_{|\lambda|\leq1}\) in
\eqref{eq:disc-modulus} and the use of a boundary average in
\eqref{eq:delta} serve different purposes.  The inequality
\eqref{eq:modulus-comparison} is always valid, but equality or an
equivalence between \(\UC(X)<2\) and complex uniform convexity does
not follow.  \rev{Such conclusions require a separate estimate
relating the disc maximum to the circular mean.}
\end{remark}

\section{Stability and Banach--Mazur estimates}

We record three stability facts.  All Banach spaces in this section
are assumed to be nonzero.

\begin{proposition}
\label{prop:stability}
\begin{enumerate}[label=\textup{(\roman*)}]
\item If \(Y\) is a nonzero closed subspace of \(X\), then
\[
 \UC(Y)\leq\UC(X),\qquad
 \LC(Y)\geq\LC(X),\qquad
 \BC(Y)\leq\BC(X).
\]
\item For a nonempty family \((X_j)_{j\in J}\) of nonzero
complex Banach spaces,
\[
 \UC\left(\bigoplus_{j\in J}^{\ell_2}X_j\right)
 =\sup_j\UC(X_j),
\]
\[
 \LC\left(\bigoplus_{j\in J}^{\ell_2}X_j\right)
 =\inf_j\LC(X_j).
\]
\item If \(T:X\to Y\) is a complex-linear isomorphism and
\(d_T=\|T\|\|T^{-1}\|\), then
\[
 d_T^{-2}\UC(X)\leq\UC(Y)\leq d_T^2\UC(X),
\]
\[
 d_T^{-2}\LC(X)\leq\LC(Y)\leq d_T^2\LC(X).
\]
\end{enumerate}
\end{proposition}

\begin{proof}
For (i), the inherited norm gives \(M_Y(x,y)=M_X(x,y)\) for
\(x,y\in Y\).  Thus the quotient defining the circular constants
is unchanged, while its domain is restricted from nonzero pairs in
\(X\times X\) to those in \(Y\times Y\).  Taking the supremum
and infimum yields
\(\UC(Y)\leq\UC(X)\) and \(\LC(Y)\geq\LC(X)\).
Since the lower constants are positive by Proposition~\ref{prop:range},
\[
 \BC(Y)=\max\left\{\UC(Y),\frac1{\LC(Y)}\right\}
 \leq\max\left\{\UC(X),\frac1{\LC(X)}\right\}=\BC(X).
\]

For (ii), put \(Z=\bigoplus_{j\in J}^{\ell_2}X_j\), and take
\(x=(x_j),y=(y_j)\in Z\).  By square summability, for each
integer \(k\geq1\), only finitely many indices satisfy
\(\|x_j\|^2+\|y_j\|^2\geq1/k\).  Their joint coordinate support
is therefore at most countable.
The coordinate squared norms are continuous and nonnegative on
\(\T\), so Tonelli's theorem on this countable support gives
\[
 \begin{aligned}
 M_Z(x,y)^2
 &=\int_{\T}\sum_{j\in J}\|x_j+\zeta y_j\|^2\dm(\zeta)\\
 &=\sum_{j\in J}M_{X_j}(x_j,y_j)^2.
 \end{aligned}
\]
The integrand is bounded by
\(2(\|x\|_Z^2+\|y\|_Z^2)\), so both sides are finite.
Also,
\[
 \|x\|_Z^2+\|y\|_Z^2
 =\sum_{j\in J}\bigl(\|x_j\|^2+\|y_j\|^2\bigr).
\]
Set \(U=\sup_{j\in J}\UC(X_j)\) and
\(L=\inf_{j\in J}\LC(X_j)\).  Proposition~\ref{prop:range}
ensures \(\tfrac12\leq L\leq1\leq U\leq2\).
Summing the defining coordinate bounds, which also hold when
\(x_j=y_j=0\), yields
\[
 L\bigl(\|x\|_Z^2+\|y\|_Z^2\bigr)
 \leq M_Z(x,y)^2
 \leq U\bigl(\|x\|_Z^2+\|y\|_Z^2\bigr).
\]
For nonzero pairs, division by the positive denominator followed
by the appropriate infimum or supremum gives
\(\LC(Z)\geq L\) and \(\UC(Z)\leq U\).
Each \(X_j\) is isometric to a closed coordinate subspace of
\(Z\), so (i) gives \(\LC(Z)\leq\LC(X_j)\) and
\(\UC(Z)\geq\UC(X_j)\) for every \(j\).
Taking the infimum and supremum over \(j\), respectively, proves
the reverse inequalities and hence both identities.

For (iii), the operator norm estimates imply, for every \(z\in X\),
\[
 \|T^{-1}\|^{-1}\|z\|_X
 \leq\|Tz\|_Y\leq\|T\|\,\|z\|_X.
\]
Complex linearity gives \(T(x+\zeta y)=Tx+\zeta Ty\).
Applying these norm estimates inside the circular mean and to each
coordinate of a pair therefore yields
\[
 \begin{aligned}
 \|T^{-1}\|^{-2}M_X(x,y)^2
 &\leq M_Y(Tx,Ty)^2\leq\|T\|^2 M_X(x,y)^2,\\
 \|T^{-1}\|^{-2}\bigl(\|x\|_X^2+\|y\|_X^2\bigr)
 &\leq\|Tx\|_Y^2+\|Ty\|_Y^2\\
 &\leq\|T\|^2\bigl(\|x\|_X^2+\|y\|_X^2\bigr).
 \end{aligned}
\]
Consequently, for every \((x,y)\ne(0,0)\),
\[
 d_T^{-2}\frac{M_X(x,y)^2}{\|x\|_X^2+\|y\|_X^2}
 \leq\frac{M_Y(Tx,Ty)^2}{\|Tx\|_Y^2+\|Ty\|_Y^2}
 \leq d_T^2\frac{M_X(x,y)^2}{\|x\|_X^2+\|y\|_X^2}.
\]
The map \((x,y)\mapsto(Tx,Ty)\) is a bijection between the
nonzero pair domains.  Taking the supremum throughout yields the
two bounds for \(\UC(Y)\); taking the infimum throughout yields
the two bounds for \(\LC(Y)\).
\end{proof}

\begin{corollary}
\label{cor:BM}
Let \(X\) be an \(n\)-dimensional complex Banach space, with
\(n\geq1\).  Then
\[
 \BC(X)\leq d(X,\ell_2^n)^2,
 \qquad
 d(X,\ell_2^n)\geq\sqrt{\BC(X)},
 \tag{6.1}\label{eq:BM}
\]
where \(d\) is the Banach--Mazur distance, defined for
finite-dimensional complex Banach spaces \(X\) and \(Y\) of the
same dimension by
\[
 d(X,Y)=\inf\bigl\{\|T\|\,\|T^{-1}\|:
 T:X\to Y\text{ is a complex-linear isomorphism}\bigr\},
\]
with \(\|T\|\) and \(\|T^{-1}\|\) denoting the corresponding
operator norms.  Thus \(d(X,Y)\) is the infimum of the
multiplicative distortions of all such isomorphisms.
\end{corollary}

\begin{proof}
Put \(H=\ell_2^n\).  Since \(H\) is a complex Hilbert space,
Proposition~\ref{prop:range} gives \(\UC(H)=\LC(H)=1\).
The spaces \(X\) and \(H\) have the same finite complex dimension,
so the set of complex-linear isomorphisms \(T:X\to H\) is nonempty.
For each such map, put \(d_T=\|T\|\|T^{-1}\|\); then
\(d_T\geq\|I_X\|=1\).
Proposition~\ref{prop:stability}(iii) gives
\[
 \UC(X)\leq d_T^2\UC(H)=d_T^2,
 \qquad
 \LC(X)\geq d_T^{-2}\LC(H)=d_T^{-2}.
\]
The positivity of the lower constant therefore implies
\[
 \BC(X)=\max\left\{\UC(X),\frac1{\LC(X)}\right\}
 \leq d_T^2.
\]
Since this holds for every such \(T\), taking the infimum gives
\[
 \BC(X)\leq\inf_T d_T^2
 =\left(\inf_T d_T\right)^2=d(X,H)^2.
\]
Continuity and monotonicity of squaring on \([0,\infty)\)
justify the equality; no minimizing isomorphism is required.
The second bound in~\eqref{eq:BM} follows by taking nonnegative
square roots.
\end{proof}

\section{Computations in classical spaces}

The examples in this section illustrate three different mechanisms:
disjoint supports, many available phases, and comparison with the
classical constant.

\subsection{\texorpdfstring{\(\ell_1\), \(\ell_\infty\) and \(C(K)\)}
{ell1, ell-infinity and C(K)}}

\begin{proposition}\label{prop:l1}
If \(X\) contains an isometric copy of complex \(\ell_1^2\), then
\[
 \UC(X)=2.
 \tag{7.1}\label{eq:l1}
\]
In particular, \(\UC(\ell_1^n)=2\) for \(n\geq2\), and
\(\UC(\ell_1)=2\).
\end{proposition}

\begin{proof}
Use the coordinate pair from Example \ref{ex:endpoints}(i) and the
universal upper bound.
\end{proof}

For finite-dimensional \(\ell_\infty\), the upper constant records
the number of available coordinate phases.  We use one elementary
planar fact.

\begin{lemma}\label{lem:polygon}
Let \(n\geq2\) be an integer, let \(u_1,\ldots,u_n\in\T\), and put
\[
 h(\zeta)=\max\{0,\operatorname{Re}(\zeta\overline{u_1}),\ldots,
 \operatorname{Re}(\zeta\overline{u_n})\}.
\]
Then
\[
 \int_{\T}h(\zeta)\dm(\zeta)
 \leq\frac{n}{\pi}\sin\frac{\pi}{n}.
 \tag{7.2}\label{eq:polygon}
\]
Equality is attained when the \(u_j\)'s are equally spaced.
\end{lemma}

\begin{proof}
Write \(u_j=e^{i\theta_j}\), ordering the arguments so that
\(\theta_1\leq\cdots\leq\theta_n<\theta_1+2\pi\), and put
\(\theta_{n+1}=\theta_1+2\pi\) and
\(\alpha_j=\theta_{j+1}-\theta_j\).  Then
\(\alpha_j\geq0\) and \(\sum_{j=1}^n\alpha_j=2\pi\);
repeated phase points are allowed.
For \(0\leq s\leq\alpha_j\), a shortest circular path from
\(e^{i(\theta_j+s)}\) to a phase point must meet an endpoint of
this gap.  Its distance to the nearest phase point is therefore
\(\min\{s,\alpha_j-s\}\leq\pi\).  Since cosine decreases on
\([0,\pi]\),
\[
 h(e^{i(\theta_j+s)})
 =\max\{0,\cos(\min\{s,\alpha_j-s\})\}.
\]
Thus, including gaps of length zero or greater than \(\pi\),
\[
 \int_{\theta_j}^{\theta_{j+1}}h(e^{i\theta})\,d\theta
 =2\int_0^{\alpha_j/2}\max\{0,\cos s\}\,ds
 =2F(\alpha_j),
\]
where
\[
 F(a)=\sin\bigl(\min\{a/2,\pi/2\}\bigr),
 \qquad 0\leq a\leq2\pi.
\]
The derivative of \(F\) is \(\tfrac12\cos(a/2)\) on
\((0,\pi)\) and zero on \((\pi,2\pi)\); the one-sided derivatives
at \(\pi\) agree.  Hence the derivative is nonincreasing and
\(F\) is concave on \([0,2\pi]\).
Using periodicity and \(dm(e^{i\theta})=d\theta/(2\pi)\),
summing over the arcs and applying Jensen's inequality yields
\[
 \int_{\T}h(\zeta)\dm(\zeta)
 =\frac1\pi\sum_{j=1}^n F(\alpha_j)
 \leq\frac n\pi F\left(\frac{2\pi}{n}\right)
 =\frac n\pi\sin\frac\pi n,
\]
where the last equality uses \(n\geq2\).
For equally spaced phase points, all \(\alpha_j=2\pi/n\),
and equality holds throughout.
\end{proof}

\begin{theorem}
\label{thm:linfty}
For every integer \(n\geq2\),
\[
 \boxed{\displaystyle
 \UC(\ell_\infty^n)
 =1+\frac{n}{\pi}\sin\frac{\pi}{n}.}
 \tag{7.3}\label{eq:linfty}
\]
Moreover,
\[
 \LC(\ell_\infty^n)=\frac12,\qquad
 \BC(\ell_\infty^n)=2.
 \tag{7.4}\label{eq:linfty-lower}
\]
\end{theorem}

\begin{proof}
Take \((x,y)\ne(0,0)\) in \(\ell_\infty^n\times\ell_\infty^n\)
and put
\[
 A=\|x\|_\infty,\qquad B=\|y\|_\infty,\qquad D=A^2+B^2>0.
\]
For \(x_jy_j\ne0\), choose \(u_j\in\T\) such that
\(\overline{x_j}y_j=|x_j||y_j|\overline{u_j}\);
if \(x_jy_j=0\), choose \(u_j\) arbitrarily.
For every coordinate and every \(\zeta\in\T\),
\[
 |x_j+\zeta y_j|^2
 =|x_j|^2+|y_j|^2
   +2|x_j||y_j|\operatorname{Re}(\zeta\overline{u_j}).
\]
We have \(|x_j|^2+|y_j|^2\leq D\) and
\(0\leq2|x_j||y_j|\leq2AB\leq D\).
Replacing the real part by its nonnegative part before using this
coefficient bound gives
\[
 \begin{aligned}
 |x_j+\zeta y_j|^2
 &\leq D+2|x_j||y_j|
                \max\{0,\operatorname{Re}(\zeta\overline{u_j})\}\\
 &\leq D\bigl(1+\max\{0,
                       \operatorname{Re}(\zeta\overline{u_j})\}\bigr).
 \end{aligned}
\]
Consequently, with \(h\) as in Lemma~\ref{lem:polygon},
\(\|x+\zeta y\|_\infty^2\leq D(1+h(\zeta))\).
Integration and \eqref{eq:polygon} yield
\[
 \frac{M_{\ell_\infty^n}(x,y)^2}{D}
 \leq1+\frac n\pi\sin\frac\pi n.
\]
Taking the supremum proves the upper bound in \eqref{eq:linfty}.

For the matching lower bound, put \(\omega_n=e^{2\pi i/n}\) and take
\[
 x=(1,\ldots,1),\qquad
 y=(1,\omega_n^{-1},\ldots,\omega_n^{-(n-1)}).
\]
Both vectors have norm one, and
\[
 \|x+e^{i\theta}y\|_\infty^2
 =2+2\max_{0\leq j<n}
       \cos\left(\theta-\frac{2\pi j}{n}\right).
\]
On the arc centered at \(2\pi j/n\) with half-length \(\pi/n\),
read modulo \(2\pi\), the maximum is the corresponding cosine.
These arcs cover the circle with disjoint interiors; ties at their
endpoints do not affect the integral.  Thus
\[
 \frac1{2\pi}\int_0^{2\pi}\max_{0\leq j<n}
       \cos\left(\theta-\frac{2\pi j}{n}\right)\,d\theta
 =\frac n{2\pi}\int_{-\pi/n}^{\pi/n}\cos s\,ds
 =\frac n\pi\sin\frac\pi n.
\]
Since \(\|x\|_\infty^2+\|y\|_\infty^2=2\), the quotient of
this pair is \(1+\frac n\pi\sin(\pi/n)\), proving \eqref{eq:linfty}.

Finally, \(n\geq2\) provides distinct coordinate vectors
\(e_1,e_2\) satisfying \(\|e_1+\zeta e_2\|_\infty=1\)
for every \(\zeta\in\T\).  Their circular quotient is \(1/2\).
The lower bound in Proposition~\ref{prop:range} therefore gives
\(\LC(\ell_\infty^n)=1/2\).  Since \(\UC(\ell_\infty^n)\leq2\),
the definition of the balanced constant gives
\(\BC(\ell_\infty^n)=2\), proving \eqref{eq:linfty-lower}.
\end{proof}

\begin{samepage}
\begin{corollary}\label{cor:infinite-square}
\[
 \UC(c_0)=\UC(\ell_\infty)=2,
 \qquad
 \LC(c_0)=\LC(\ell_\infty)=\frac12.
\]
If \(K\) is a nonempty compact Hausdorff space, then
\[
 \UC(C(K))=
 \begin{cases}
 1, & |K|=1,\\[2mm]
 1+\dfrac{|K|}{\pi}\sin\dfrac{\pi}{|K|},
     & 2\leq |K|<\infty,\\[3mm]
 2, & |K|=\infty.
 \end{cases}
 \tag{7.5}\label{eq:CK}
\]
If \(|K|\geq2\), then \(\LC(C(K))=1/2\) and
\(\BC(C(K))=2\).
\end{corollary}
\end{samepage}

\begin{proof}
For every integer \(n\geq2\) and \(X=c_0\) or \(\ell_\infty\),
the zero-extension embedding of
\(\ell_\infty^n\) and Proposition~\ref{prop:stability}(i) give
\[
 \UC(X)\geq\UC(\ell_\infty^n)
 =1+\frac n\pi\sin\frac\pi n\longrightarrow2.
\]
The universal upper bound yields \(\UC(X)=2\).  The coordinate
vectors \(e_1,e_2\) have circular quotient \(1/2\); hence
Proposition~\ref{prop:range} gives \(\LC(X)=1/2\).

If \(K\) is finite, then \(C(K)\) is isometric to
\(\ell_\infty^{|K|}\); the stated values follow from
Proposition~\ref{prop:range} when \(|K|=1\), and from
Theorem~\ref{thm:linfty} when \(|K|\geq2\).

Suppose \(K\) is infinite.  For \(n\geq2\), choose distinct
\(t_0,\ldots,t_{n-1}\in K\).  \rev{By normality, choose pairwise
disjoint open neighborhoods \(U_j\) of these points.  Urysohn's
lemma gives \(\varphi_j\in C(K)\), with
\(0\leq\varphi_j\leq1\), \(\varphi_j(t_j)=1\), and
\(\operatorname{supp}\varphi_j\subset U_j\).}
Put \(\omega_n=e^{2\pi i/n}\) and
\(V=\sum_{j=0}^{n-1}\omega_n^{-j}\varphi_j\).
Then \(\|1\|_\infty=\|V\|_\infty=1\) and, for \(\zeta\in\T\),
\[
 \|1+\zeta V\|_\infty^2
 \geq\max_{0\leq j<n}|1+\zeta\omega_n^{-j}|^2.
\]
The integral calculated in the proof of Theorem~\ref{thm:linfty}
therefore gives
\[
 \UC(C(K))\geq\frac{M_{C(K)}(1,V)^2}{2}
 \geq1+\frac n\pi\sin\frac\pi n\longrightarrow2.
\]
The universal upper bound proves \(\UC(C(K))=2\).
When \(|K|\geq2\), two disjointly supported unit functions
\(f,g\in C(K)\) satisfy \(\|f+\zeta g\|_\infty=1\) for every
\(\zeta\in\T\).  Thus their circular quotient is \(1/2\), and
Proposition~\ref{prop:range} yields \(\LC(C(K))=1/2\).
Since \(\UC(C(K))\leq2\), the definition gives \(\BC(C(K))=2\).
\end{proof}

\begin{corollary}
\label{cor:l1-lower}
For \(n\geq2\),
\[
 \LC(\ell_1^n)\geq
 \left(1+\frac n\pi\sin\frac\pi n\right)^{-1}.
 \tag{7.5a}\label{eq:l1-lower}
\]
Moreover, the Hadamard pair supported on two coordinates gives
\[
 \LC(\ell_1^n)\leq\frac12+\frac1\pi.
 \tag{7.5b}\label{eq:l1-upper}
\]
\end{corollary}

\begin{proof}
Since \((\ell_1^n)^*=\ell_\infty^n\), the duality inequality
\eqref{eq:dual} and Theorem~\ref{thm:linfty} yield
\[
 \LC(\ell_1^n)\geq\frac1{\UC(\ell_\infty^n)}
 =\left(1+\frac n\pi\sin\frac\pi n\right)^{-1}.
\]
For the upper bound, choose \(x=(1,1,0,\ldots,0)\) and
\(y=(1,-1,0,\ldots,0)\).  Then \(\|x\|_1=\|y\|_1=2\), so the
denominator of their circular quotient is \(8\).  For \(t\in[0,2\pi]\),
\[
 \begin{aligned}
 \|x+e^{it}y\|_1^2
 &=\bigl(|1+e^{it}|+|1-e^{it}|\bigr)^2\\
 &=4\left(\left|\cos\frac t2\right|
          +\left|\sin\frac t2\right|\right)^2
 =4(1+|\sin t|).
 \end{aligned}
\]
Using \(\int_0^{2\pi}|\sin t|\,dt=4\), we obtain
\[
 M_{\ell_1^n}(x,y)^2
 =\frac1{2\pi}\int_0^{2\pi}4(1+|\sin t|)\,dt
 =4+\frac8\pi.
\]
Consequently this pair has quotient \((4+8/\pi)/8=1/2+1/\pi\).
The infimum defining \(\LC(\ell_1^n)\) is no larger, proving
\eqref{eq:l1-upper}.
\end{proof}

\subsection{\texorpdfstring{\(L_p\) and \(\ell_p\)}{Lp and ell-p}}

We use the classical formula
\[
 \CNJ(L_p)=\CNJ(\ell_p)=2^{|2/p-1|}
 \qquad(1\leq p\leq\infty)
 \tag{7.6}\label{eq:classical-Lp}
\]
whenever two disjoint unit vectors are available; see
\cite{Clarkson1937} and
\rev{\cite[Corollary~4]{KatoTakahashiLebesgue1998}}.

\begin{theorem}\label{thm:Lp-low}
Let \(1\leq p\leq2\).  If \(L_p(\mu)\) contains two disjoint unit
vectors, then
\[
 \boxed{\displaystyle
 \UC(L_p(\mu))=2^{2/p-1}.}
 \tag{7.7}\label{eq:Lp-low}
\]
In particular, for every \(n\geq2\),
\[
 \UC(\ell_p^n)=\UC(\ell_p)=2^{2/p-1}.
 \tag{7.8}\label{eq:lp-low}
\]
\end{theorem}

\begin{proof}
Let \(u,v\in L_p(\mu)\) be disjoint unit vectors.  Since their
supports are disjoint up to null sets, for every \(\zeta\in\T\),
\[
 \|u+\zeta v\|_p^p
 =\int(|u|^p+|v|^p)\,d\mu=2.
\]
Hence \(M_{L_p}(u,v)^2=2^{2/p}\), while
\(\|u\|_p^2+\|v\|_p^2=2\).  The defining supremum therefore gives
\(\UC(L_p(\mu))\geq2^{2/p-1}\).  Conversely,
\eqref{eq:comparison} and \eqref{eq:classical-Lp} imply
\[
 \UC(L_p(\mu))\leq\CNJ(L_p(\mu))=2^{2/p-1}.
\]
The two bounds prove \eqref{eq:Lp-low}.  In \(\ell_p^n\),
\(n\geq2\), and in \(\ell_p\), two distinct coordinate vectors
provide the required disjoint unit vectors.
\end{proof}

The other side of the \(L_p\)-scale is also exact.

\begin{theorem}
\label{thm:Lp-lower-high}
Let \(2\leq p\leq\infty\).  If \(L_p(\mu)\) contains two disjoint
unit vectors, then
\[
 \boxed{\displaystyle
 \LC(L_p(\mu))=2^{2/p-1},}
 \tag{7.9}\label{eq:Lp-lower-high}
\]
with the convention \(2^{2/\infty-1}=1/2\).  The same formula holds
for \(\ell_p^n\), \(n\geq2\), and \(\ell_p\).
\end{theorem}

\begin{proof}
Let \(u,v\in L_p(\mu)\) be disjoint unit vectors.  For
\(2\leq p<\infty\), disjointness gives
\(\|u+\zeta v\|_p^p=2\) for every \(\zeta\in\T\).  When
\(p=\infty\), their disjoint supports give
\(\|u+\zeta v\|_\infty=\max\{\|u\|_\infty,\|v\|_\infty\}=1\).
Thus, with \(2^{1/\infty}=1\),
\[
 \frac{M_{L_p}(u,v)^2}{\|u\|_p^2+\|v\|_p^2}
 =\frac{2^{2/p}}2=2^{2/p-1}.
\]
Taking the infimum yields \(\LC(L_p(\mu))\leq2^{2/p-1}\).  On the
other hand, \eqref{eq:comparison} and \eqref{eq:classical-Lp} give
\[
 \LC(L_p(\mu))\geq\frac1{\CNJ(L_p(\mu))}
 =2^{2/p-1}.
\]
This proves \eqref{eq:Lp-lower-high}; the sequence-space statements
follow by taking two distinct coordinate vectors.
\end{proof}

For the two remaining one-sided problems, a rotated Hadamard pair
gives a useful explicit bound.  For \(1\leq p<\infty\), put
\[
 B_p=
 2^{1-2/p}\frac2\pi
 \int_0^{\pi/2}
 \bigl(\cos^p s+\sin^p s\bigr)^{2/p}\,ds.
 \tag{7.10}\label{eq:Bp}
\]

\begin{theorem}
\label{thm:Lp-high}
If \(1\leq p\leq2\) and two disjoint unit vectors are available,
then
\[
 2^{1-2/p}\leq\LC(L_p(\mu))\leq B_p.
 \tag{7.11}\label{eq:Lp-low-lower}
\]
If \(2\leq p<\infty\) and \(L_p(\mu)\) contains two disjoint unit
vectors, then
\[
 B_p\leq\UC(L_p(\mu))\leq2^{1-2/p}.
 \tag{7.12}\label{eq:Lp-high}
\]
The same bounds hold for \(\ell_p^n\), \(n\geq2\), and
\(\ell_p\).  Moreover,
\[
 B_2=1,\qquad
 \lim_{p\to\infty}B_p=1+\frac2\pi.
 \tag{7.13}\label{eq:Bp-limit}
\]
\end{theorem}

\begin{proof}
Let \(u,v\) be disjoint unit vectors, and set \(x=u+v\) and
\(y=u-v\).  For \(1\leq p<\infty\), disjointness implies
\(\|x\|_p=\|y\|_p=2^{1/p}\), and hence
\(\|x\|_p^2+\|y\|_p^2=2^{1+2/p}\).  Moreover,
\[
 \|x+e^{it}y\|_p^p=|1+e^{it}|^p+|1-e^{it}|^p.
\]
Using \(|1+e^{it}|=2|\cos(t/2)|\) and
\(|1-e^{it}|=2|\sin(t/2)|\), we obtain
\[
 \|x+e^{it}y\|_p^2
 =4\bigl(|\cos(t/2)|^p+|\sin(t/2)|^p\bigr)^{2/p}.
\]
The change of variables \(s=t/2\) and symmetry about \(\pi/2\)
give the circular mean
\[
 \begin{aligned}
 M_{L_p}(x,y)^2
 &=\frac4\pi\int_0^\pi
    \bigl(|\cos s|^p+|\sin s|^p\bigr)^{2/p}\,ds\\
 &=\frac8\pi\int_0^{\pi/2}
    (\cos^p s+\sin^p s)^{2/p}\,ds.
 \end{aligned}
\]
Dividing by \(2^{1+2/p}\), we find
\[
 \frac{M_{L_p}(x,y)^2}{\|x\|_p^2+\|y\|_p^2}
 =2^{1-2/p}\frac2\pi\int_0^{\pi/2}
   (\cos^p s+\sin^p s)^{2/p}\,ds=B_p.
\]
For \(1\leq p\leq2\), the infimum definition and
\eqref{eq:comparison} therefore give
\[
 2^{1-2/p}=\frac1{\CNJ(L_p(\mu))}
 \leq\LC(L_p(\mu))\leq B_p.
\]
For \(2\leq p<\infty\), the supremum definition and the same
comparison give
\[
 B_p\leq\UC(L_p(\mu))\leq\CNJ(L_p(\mu))=2^{1-2/p}.
\]
The same construction with two coordinate vectors proves the
sequence-space bounds.

At \(p=2\), the integrand in \eqref{eq:Bp} equals \(1\), giving
\(B_2=1\).  For \(p\geq2\) and \(s\in[0,\pi/2]\),
\(\cos^p s+\sin^p s\leq\cos^2s+\sin^2s=1\).  Thus the integrand
is bounded by \(1\) and converges pointwise to
\(\max\{\cos^2s,\sin^2s\}\).  Dominated convergence
and \(2^{1-2/p}\to2\) yield
\[
 \begin{aligned}
 \lim_{p\to\infty}B_p
 &=\frac4\pi\left(\int_0^{\pi/4}\cos^2s\,ds
       +\int_{\pi/4}^{\pi/2}\sin^2s\,ds\right)\\
 &=\frac4\pi\left(\frac\pi4+\frac12\right)=1+\frac2\pi.
 \end{aligned}
\]
\end{proof}

\begin{remark}
It is natural to ask whether \(B_p\) is the exact value of
\(\UC(\ell_p^2)\) for \(p>2\), and of \(\LC(\ell_p^2)\) for \rev{\(1\leq p<2\)}.
No such equality is used here; \eqref{eq:Lp-low-lower} and
\eqref{eq:Lp-high} are the established bounds.
\end{remark}

\subsection{Schatten classes and operator spaces}

Let \(S_p(H)\) denote the Schatten class on a \rev{complex Hilbert
space} \(H\), with \(S_\infty(H)=K(H)\) equipped with the operator norm.
The von Neumann--Jordan constants of Schatten classes agree with the
corresponding \(\ell_p\) values when \(\dim H\geq2\); see
\rev{\cite[Corollary~4]{KatoTakahashiLebesgue1998}}.  Orthogonal rank-one projections
play the role of disjoint coordinate vectors.

\begin{theorem}
\label{thm:Schatten}
Let \(\dim H\geq2\).  Then
\[
 \UC(S_p(H))=2^{2/p-1}\qquad(1\leq p\leq2),
 \tag{7.14}\label{eq:Schatten-upper}
\]
and
\[
 \LC(S_p(H))=2^{2/p-1}\qquad(2\leq p\leq\infty).
 \tag{7.15}\label{eq:Schatten-lower}
\]
For the unresolved sides one has the same Hadamard bounds as in
\eqref{eq:Lp-low-lower} and \eqref{eq:Lp-high}.
\rev{At \(p=\infty\), the corresponding upper-side bound is
\(1+2/\pi\leq\UC(S_\infty(H))\leq2\).}
\end{theorem}

\begin{proof}
Choose orthogonal rank-one projections \(P,Q\) on \(H\).  For every
\(\zeta\in\T\), the nonzero singular values of \(P+\zeta Q\) are
\(1,1\).  Thus \(\|P\|_p=\|Q\|_p=1\) and
\(\|P+\zeta Q\|_p=2^{1/p}\), with \(2^{1/\infty}=1\).  In particular,
\[
 \frac{M_{S_p(H)}(P,Q)^2}{\|P\|_p^2+\|Q\|_p^2}
 =\frac{2^{2/p}}2=2^{2/p-1}.
\]
For \(1\leq p\leq2\), this gives the lower bound for \(\UC\),
while \eqref{eq:comparison} and the classical Schatten formula give
\(\UC(S_p(H))\leq\CNJ(S_p(H))=2^{2/p-1}\).  For
\(2\leq p\leq\infty\), the same pair gives the upper bound for
\(\LC\), and
\[
 \LC(S_p(H))\geq\frac1{\CNJ(S_p(H))}
 =2^{2/p-1}.
\]
This proves \eqref{eq:Schatten-upper} and \eqref{eq:Schatten-lower}.

For the remaining bounds, put \(A=P+Q\) and \(B=P-Q\).  On
\(E=\operatorname{ran}P\oplus\operatorname{ran}Q\), the singular
values of \(A+\zeta B\), counting zeros, are \(|1+\zeta|\) and
\(|1-\zeta|\); the operator vanishes on \(E^\perp\).  In particular,
one of these values is zero when \(\zeta=\pm1\).  Also,
\(\|A\|_p=\|B\|_p=2^{1/p}\).
For \(1\leq p<\infty\), the calculation in Theorem~\ref{thm:Lp-high}
shows that their circular quotient is \(B_p\).  Hence, using
\eqref{eq:comparison} and the classical Schatten formula,
\[
 \begin{aligned}
 2^{1-2/p}&\leq\LC(S_p(H))\leq B_p &&(1\leq p\leq2),\\
 B_p&\leq\UC(S_p(H))\leq2^{1-2/p} &&(2\leq p<\infty).
 \end{aligned}
\]
For \(p=\infty\), putting \(\zeta=e^{it}\) gives
\(\max\{|1+\zeta|^2,|1-\zeta|^2\}=2+2|\cos t|\).  The quotient
of this pair is therefore
\[
 \frac12\int_{\T}\max\{|1+\zeta|^2,|1-\zeta|^2\}\dm(\zeta)
 =1+\frac2\pi=\lim_{p\to\infty}B_p,
\]
and the universal upper bound is \(2\).
\end{proof}

\begin{corollary}
\label{cor:operators}
If \(H\) is an \rev{infinite-dimensional complex Hilbert space}, then
\[
 \UC(K(H))=\UC(B(H))=2,\qquad
 \LC(K(H))=\LC(B(H))=\frac12.
 \tag{7.16}\label{eq:operator-extremes}
\]
For \(n\geq2\), the matrix algebra \(M_n\) with the operator norm
satisfies
\[
 1+\frac n\pi\sin\frac\pi n
 \leq\UC(M_n)\leq2,\qquad \LC(M_n)=\frac12.
 \tag{7.17}\label{eq:matrix-bounds}
\]
\end{corollary}

\begin{proof}
Choose an orthonormal sequence \((e_j)\) in \(H\), and let
\(E=\overline{\operatorname{span}}\{e_j:j\geq1\}\).
For \(a=(a_j)\in\ell_\infty\), define the diagonal operator
\(D_a e_j=a_j e_j\) on \(E\) and \(D_a=0\) on \(E^\perp\).
Then \(\|D_a\|=\|a\|_\infty\); if \(a\in c_0\), its finite-rank
truncations converge in operator norm, so \(D_a\) is compact.
Thus \(K(H)\) and \(B(H)\) contain closed isometric copies of
\(c_0\) and \(\ell_\infty\), respectively.  Corollary~\ref{cor:infinite-square}
and Proposition~\ref{prop:stability}(i) give, for either operator space \(X\),
\[
 \UC(X)\geq2,\qquad \LC(X)\leq\frac12.
\]
The reverse inequalities follow from Proposition~\ref{prop:range},
proving \eqref{eq:operator-extremes}.

The diagonal subspace of \(M_n\) is isometric to \(\ell_\infty^n\).
Hence Theorem~\ref{thm:linfty} and subspace monotonicity yield
\[
 1+\frac n\pi\sin\frac\pi n=\UC(\ell_\infty^n)
 \leq\UC(M_n)\leq2.
\]
Moreover, the diagonal matrix units \(E_{11},E_{22}\) have norm one
and satisfy \(\|E_{11}+\zeta E_{22}\|=1\) for every \(\zeta\in\T\).
Consequently,
\[
 \frac{M_{M_n}(E_{11},E_{22})^2}
 {\|E_{11}\|^2+\|E_{22}\|^2}
 =\frac{\int_\T 1\,dm(\zeta)}{2}=\frac12.
\]
Together with the universal lower bound, this gives
\(\LC(M_n)=1/2\) and completes \eqref{eq:matrix-bounds}.
\end{proof}

\begin{remark}
The finite-dimensional boundary behavior is already extremal:
\(\BC(\ell_\infty^n)=2\) for every \(n\geq2\), while
\[
 \UC(\ell_\infty^n)=1+\frac n\pi\sin\frac\pi n\uparrow2.
\]
Thus upper circular constants can be arbitrarily close to \(2\) in
finite dimensions.  The equalities for \(c_0\) and \(\ell_\infty\)
in Corollary \ref{cor:infinite-square} show that this boundary value is
stable in two standard non-reflexive spaces.
\end{remark}

\section{Uniform nonsquareness and fixed points}

A Banach space \(X\) is uniformly nonsquare if there is
\(\delta>0\) such that
\[
 \min\{\|x+y\|,\|x-y\|\}\leq2(1-\delta)
 \qquad(x,y\in S_X).
 \tag{8.1}\label{eq:UNS}
\]
Equivalently, its James constant
\[
 J(X)=\sup_{x,y\in S_X}
 \min\{\|x+y\|,\|x-y\|\}
 \tag{8.2}\label{eq:James}
\]
satisfies \(J(X)<2\).
The relation between von Neumann--Jordan constants, uniform
nonsquareness, and normal structure is developed in
\cite{KatoTakahashi1998,KatoMaligrandaTakahashi2001,Saejung2006}.

\begin{theorem}
\label{thm:James-quantitative}
Every complex Banach space satisfies
\[
 \boxed{\displaystyle
 \UC(X)\geq1+\frac2\pi
 -4\sqrt{2\bigl(2-J(X)\bigr)}.}
 \tag{8.3}\label{eq:James-lower}
\]
Consequently, if
\(\gamma=1+2/\pi-\UC(X)>0\), then
\[
 J(X)\leq2-\frac{\gamma^2}{32}.
 \tag{8.4}\label{eq:James-gap}
\]
Thus the nonsquareness parameter in \eqref{eq:UNS} may be chosen at
least \(\gamma^2/64\).
\end{theorem}

\begin{proof}
Fix \(x,y\in S_X\), and set
\[
 r=\min\{\|x+y\|,\|x-y\|\},\qquad \eta=2-r\geq0.
\]
By Hahn--Banach, choose \(f,g\in S_{X^*}\) so that
\(f(x+y)=\|x+y\|\) and \(g(x-y)=\|x-y\|\), with real nonnegative
values.  If one of \(x+y,x-y\) is zero, choose the corresponding
unit functional arbitrarily.

Since \(|f(x)|,|f(y)|\leq1\) and
\(\operatorname{Re}f(x)+\operatorname{Re}f(y)\geq2-\eta\),
both real parts are at least \(1-\eta\).  For \(|z|\leq1\),
\( |z-1|^2\leq2(1-\operatorname{Re}z)\); hence
\[
 |f(x)-1|,|f(y)-1|\leq\sqrt{2\eta}.
 \tag{8.5}\label{eq:f-error}
\]
Likewise, \(\operatorname{Re}g(x)-\operatorname{Re}g(y)\geq2-\eta\)
and the two summands on the left are at most \(1\).  Applying the
same estimate to \(g(x)\) and \(-g(y)\) gives
\[
 |g(x)-1|,|g(y)+1|\leq\sqrt{2\eta}.
 \tag{8.6}\label{eq:g-error}
\]

Put \(d=2\sqrt{2\eta}\) and
\(h(\zeta)=\max\{|1+\zeta|,|1-\zeta|\}\).  For \(\zeta\in\T\),
\eqref{eq:f-error}--\eqref{eq:g-error} and the reverse triangle
inequality give
\[
 \begin{aligned}
 |f(x)+\zeta f(y)|&\geq|1+\zeta|-d,\\
 |g(x)+\zeta g(y)|&\geq|1-\zeta|-d.
 \end{aligned}
\]
Because \(\|f\|=\|g\|=1\), it follows that
\(\|x+\zeta y\|\geq h(\zeta)-d\), and hence
\(\|x+\zeta y\|^2\geq(h(\zeta)-d)_+^2\), where
\(t_+=\max\{t,0\}\).  For \(a,d\geq0\), \((a-d)_+^2\geq a^2-2da\); also \(h\leq2\) and
\[
 \int_{\T}h(\zeta)^2\dm(\zeta)
 =\frac1{2\pi}\int_0^{2\pi}(2+2|\cos t|)\,dt
 =2+\frac4\pi.
\]
It follows that
\[
 M_X(x,y)^2\geq\int_{\T}h^2\,dm-2d\int_{\T}h\,dm
 \geq2+\frac4\pi-4d.
\]
Since \(\|x\|^2+\|y\|^2=2\), this yields
\[
 \frac{M_X(x,y)^2}{2}
 \geq1+\frac2\pi-2d
 =1+\frac2\pi-4\sqrt{2\eta}.
\]
Choose pairs \((x_n,y_n)\) for which \(r_n\to J(X)\).
Their quotients are bounded above by \(\UC(X)\), so passing to the
limit proves \eqref{eq:James-lower}.  If
\(\gamma=1+2/\pi-\UC(X)>0\), then
\(\gamma\leq4\sqrt{2(2-J(X))}\).  Squaring yields
\eqref{eq:James-gap}.  Indeed, every \(x,y\in S_X\) satisfies
\(\min\{\|x+y\|,\|x-y\|\}\leq J(X)\leq2(1-\gamma^2/64)\),
so \eqref{eq:UNS} holds with \(\delta=\gamma^2/64\).
\end{proof}

The circular threshold below is smaller than \(2\), because averaging
reduces the two isolated peaks at the phases \(1\) and \(-1\).

\begin{theorem}\label{thm:UNS}
If \(X\) is not uniformly nonsquare, then
\[
 \UC(X)\geq1+\frac2\pi.
 \tag{8.7}\label{eq:UNS-lower}
\]
Equivalently,
\[
 \boxed{\displaystyle
 \UC(X)<1+\frac2\pi
 \quad\Longrightarrow\quad
 X\text{ is uniformly nonsquare}.}
 \tag{8.8}\label{eq:UNS-criterion}
\]
The number \(1+2/\pi\) is optimal.
\end{theorem}

\begin{proof}
If \(X\) is not uniformly nonsquare, choose \(x_n,y_n\in S_X\)
with \(\min\{\|x_n+y_n\|,\|x_n-y_n\|\}>2-1/n\).  The triangle
inequality bounds both norms by \(2\), so
\[
 \|x_n+y_n\|\longrightarrow2,\qquad
 \|x_n-y_n\|\longrightarrow2.
 \tag{8.9}\label{eq:square-seq}
\]
Choose norm-one functionals \(f_n,g_n\) with
\(f_n(x_n+y_n)=\|x_n+y_n\|\) and
\(g_n(x_n-y_n)=\|x_n-y_n\|\).  Put
\(\eta_n=2-\min\{\|x_n+y_n\|,\|x_n-y_n\|\}\to0\).
Equations \eqref{eq:f-error} and \eqref{eq:g-error}, with
\(\eta=\eta_n\), imply
\[
 f_n(x_n)\to1,\quad f_n(y_n)\to1,\qquad
 g_n(x_n)\to1,\quad g_n(y_n)\to-1.
 \tag{8.10}\label{eq:functional-limits}
\]
Consequently, for each \(\zeta\in\T\),
\[
 \begin{aligned}
 \|x_n+\zeta y_n\|
 &\geq\max\{|f_n(x_n)+\zeta f_n(y_n)|,
               |g_n(x_n)+\zeta g_n(y_n)|\},\\
 \liminf_{n\to\infty}\|x_n+\zeta y_n\|^2
 &\geq\max\{|1+\zeta|^2,|1-\zeta|^2\}.
 \end{aligned}
\]
The integrands are continuous and nonnegative, and the circular
quotient of each pair has denominator \(2\).  Fatou's lemma gives
\[
 \begin{aligned}
 \UC(X)&\geq\frac12\liminf_{n\to\infty}
      \int_{\T}\|x_n+\zeta y_n\|^2\dm(\zeta)\\
 &\geq\frac12\int_{\T}\max\{|1+\zeta|^2,|1-\zeta|^2\}\dm(\zeta)\\
 &=\frac12\int_0^{2\pi}(2+2|\cos t|)\frac{dt}{2\pi}
  =1+\frac2\pi.
 \end{aligned}
\]
The implication \eqref{eq:UNS-criterion} is the contrapositive.
For sharpness, Theorem~\ref{thm:linfty} gives
\(\UC(\ell_\infty^2)=1+2/\pi\).  The unit vectors
\rev{\(x=(1,1)\) and \(y=(1,-1)\)} satisfy \(\|x+y\|_\infty=
\|x-y\|_\infty=2\); thus \(\ell_\infty^2\) is not uniformly
nonsquare, and the threshold cannot be increased.
\end{proof}

\begin{corollary}\label{cor:FPP}
If
\[
 \UC(X)<1+\frac2\pi,
\]
then every nonexpansive self-map of every nonempty bounded closed
convex subset of \(X\) has a fixed point.
\end{corollary}

\begin{proof}
By Theorem~\ref{thm:UNS}, the hypothesis implies that \(X\) is
uniformly nonsquare.  The fixed point theorem of García-Falset,
Llorens-Fuster and Mazcuñan-Navarro \cite{GarciaFalset2006} then
applies: every nonexpansive self-map of a nonempty bounded closed
convex subset of a uniformly nonsquare Banach space has a fixed point.
For complex \(X\), apply the theorem to the underlying real Banach
space; the norm, convex subsets, and nonexpansiveness are unchanged.
\end{proof}

\begin{remark}
For \(1<p<\infty\), the \(L_p\)-spaces already have the fixed point
property by uniform convexity.  The value of Corollary \ref{cor:FPP}
is that it applies to an arbitrary complex Banach space through one
computable circular estimate.
\end{remark}

\section{Complex interpolation and dimension-sensitive bounds}
\label{sec:interpolation}

The Fourier synthesis representation makes complex interpolation
particularly effective.  We first recall the notation used here;
see \cite[Sections~2.3 and~4.1]{BerghLofstrom1976}.
A couple \((X_0,X_1)\) is \emph{compatible} if both spaces embed
continuously into a common Hausdorff topological vector space, and
it is \emph{regular} if \(X_0\cap X_1\) is dense in each \(X_j\).
Their sum is equipped with
\[
 \|x\|_{X_0+X_1}
 =\inf_{x=x_0+x_1}
   \bigl(\|x_0\|_{X_0}+\|x_1\|_{X_1}\bigr),\qquad x_j\in X_j.
\]
Let \(\mathcal S=\{z\in\C:0\leq\operatorname{Re}z\leq1\}\), and let
\(\mathcal F(X_0,X_1)\) consist of the bounded continuous
\((X_0+X_1)\)-valued functions on \(\mathcal S\) that are holomorphic
in its interior and whose boundary maps \(t\mapsto F(j+it)\) are
continuous into \(X_j\) and tend to zero there as \(|t|\to\infty\),
for \(j=0,1\).  Set
\[
 \|F\|_{\mathcal F}
 =\max_{j=0,1}\sup_{t\in\R}\|F(j+it)\|_{X_j}.
\]
For \(0<\theta<1\), the \emph{Calder\'on complex interpolation space}
and its norm are
\[
 [X_0,X_1]_\theta=\{F(\theta):F\in\mathcal F(X_0,X_1)\},
 \qquad
 \|x\|_{[X_0,X_1]_\theta}
 =\inf_{\substack{F\in\mathcal F(X_0,X_1)\\F(\theta)=x}}
    \|F\|_{\mathcal F}.
\]
The complex interpolation theorem for linear operators
\cite[Theorem~4.1.2]{BerghLofstrom1976} states that if
\((A_0,A_1)\) and \((B_0,B_1)\) are compatible couples of complex
Banach spaces and \(T:A_0+A_1\to B_0+B_1\) is a linear operator
whose restrictions \(T:A_j\to B_j\) are bounded, with
\(M_j=\|T\|_{A_j\to B_j}\), then for every \(0<\theta<1\),
\[
 T:[A_0,A_1]_\theta\longrightarrow[B_0,B_1]_\theta,
 \qquad
 \|T\|_{[A_0,A_1]_\theta\to[B_0,B_1]_\theta}
 \leq M_0^{1-\theta}M_1^\theta.
\]
In particular, the two endpoint restrictions must arise from the same
linear operator.  We apply this theorem together with the isometric
interpolation identities for finite \(\ell_2\)-sums and Bochner
\(L_2\)-spaces.

\begin{theorem}
\label{thm:interpolation}
Let \((X_0,X_1)\) be a regular compatible couple of nonzero complex
Banach spaces, let \(0<\theta<1\), and set
\(X_\theta=[X_0,X_1]_\theta\).  Then
\[
 \boxed{\displaystyle
 \UC(X_\theta)
 \leq \UC(X_0)^{1-\theta}\UC(X_1)^\theta.}
 \tag{9.1}\label{eq:interpolation}
\]
\end{theorem}

\begin{proof}
For \(j=0,1\), let \(A_j=X_j\oplus_2X_j\) and
\(Y_j=L_2(\T;X_j)\), and define
\[
 S_j:A_j\longrightarrow Y_j,\qquad
 S_j(x,y)(\zeta)=x+\zeta y.
\]
\rev{The domain couple embeds into
\((X_0+X_1)\oplus_2(X_0+X_1)\), and the range couple embeds
continuously into \(L_2(\T;X_0+X_1)\), so both couples are
compatible.}  The two maps agree on \(A_0\cap A_1\), since both have
the same pointwise formula.  They therefore define a common linear operator
\(S:A_0+A_1\to Y_0+Y_1\).
For \((x,y)\in A_j\),
\[
 \|S_j(x,y)\|_{Y_j}^2
 =\int_{\T}\|x+\zeta y\|_{X_j}^2\,dm(\zeta)
 =M_{X_j}(x,y)^2,
 \quad \|(x,y)\|_{A_j}^2=\|x\|_{X_j}^2+\|y\|_{X_j}^2.
\]
Consequently, Proposition~\ref{prop:operator} gives
\(\|S_j\|^2=\UC(X_j)\) for \(j=0,1\).

By \cite[Theorem~5.1.2]{BerghLofstrom1976}, applied to counting
measure on two points and to \((\T,m)\), respectively, the following
identities hold with equality of norms:
\[
 [A_0,A_1]_\theta=X_\theta\oplus_2X_\theta,
 \qquad [Y_0,Y_1]_\theta=L_2(\T;X_\theta).
\]
Applying the operator interpolation theorem recalled above to
\(S:A_0+A_1\to Y_0+Y_1\) gives
\[
 \|S\|_{[A_0,A_1]_\theta\to[Y_0,Y_1]_\theta}
 \leq\|S_0\|^{1-\theta}\|S_1\|^\theta.
\]
On \((X_0\cap X_1)\times(X_0\cap X_1)\), the interpolated map
has the formula \(\Synth{X_\theta}(x,y)(\zeta)=x+\zeta y\).
This subspace is dense in \(X_\theta\oplus_2X_\theta\), by
\cite[Theorem~4.2.2]{BerghLofstrom1976}, so the same formula identifies
the map on the whole interpolation space.  Squaring and applying
Proposition~\ref{prop:operator} once more yields
\[
 \UC(X_\theta)=\|\Synth{X_\theta}\|^2
 \leq(\|S_0\|^2)^{1-\theta}(\|S_1\|^2)^\theta
 =\UC(X_0)^{1-\theta}\UC(X_1)^\theta,
\]
as required.
\end{proof}

\begin{corollary}
\label{cor:finite-p}
For each integer \(n\geq2\), put
\[
 A_n=1+\frac n\pi\sin\frac\pi n.
 \tag{9.2}\label{eq:An}
\]
With \(B_\infty:=\lim_{p\to\infty}B_p=1+2/\pi\), as in
\eqref{eq:Bp-limit}, for \(2\leq p\leq\infty\) we have
\[
 \boxed{\displaystyle
 B_p\leq\UC(\ell_p^n)\leq A_n^{1-2/p}.}
 \tag{9.3}\label{eq:finite-p}
\]
For finite \(n\) and \(2<p<\infty\), the upper bound is strictly
smaller than the dimension-free estimate \(2^{1-2/p}\).  For
\(1\leq p\leq2\), interpolation recovers the exact formula
\[
 \UC(\ell_p^n)=2^{2/p-1}.
 \tag{9.4}\label{eq:finite-low-exact}
\]
\end{corollary}

\begin{proof}
All couples used below have the common underlying vector space
\(\C^n\), so they are regular and compatible.  The Hilbert case gives
\(\UC(\ell_2^n)=1\); also, \(\UC(\ell_1^n)=2\), since
\(e_1,e_2\) give the quotient \(2\) and the universal upper bound is
\(2\).  Suppose first that \(2<p<\infty\), and set
\(\theta=1-2/p\).
The finite-dimensional interpolation identity \rev{\cite[Theorem~5.1.1]{BerghLofstrom1976}} gives, isometrically,
\[
 \ell_p^n=[\ell_2^n,\ell_\infty^n]_\theta,
 \qquad \frac1p=\frac{1-\theta}{2}.
\]
Since \(\UC(\ell_2^n)=1\) and
\(\UC(\ell_\infty^n)=A_n\) by Theorem~\ref{thm:linfty},
Theorem~\ref{thm:interpolation} yields
\[
 \UC(\ell_p^n)\leq1^{1-\theta}A_n^\theta
 =A_n^{1-2/p}.
\]
For the lower bound, take \(x=e_1+e_2\) and \(y=e_1-e_2\).
These vectors use only two coordinates, which are available since
\(n\geq2\), and satisfy
\[
 \|x\|_p^2+\|y\|_p^2=2^{1+2/p},\qquad
 \|x+e^{it}y\|_p^2
 =4\bigl(|\cos(t/2)|^p+|\sin(t/2)|^p\bigr)^{2/p}.
\]
The Hadamard calculation in Theorem~\ref{thm:Lp-high} therefore gives
\[
 \frac{M_{\ell_p^n}(x,y)^2}{\|x\|_p^2+\|y\|_p^2}
 =2^{1-2/p}\frac2\pi\int_0^{\pi/2}
   (\cos^p s+\sin^p s)^{2/p}\,ds=B_p.
\]
Taking the supremum over all nonzero pairs yields
\(B_p\leq\UC(\ell_p^n)\).
At \(p=2\), both bounds equal \(1\).  At \(p=\infty\),
Theorem~\ref{thm:linfty} gives \(\UC(\ell_\infty^n)=A_n\),
while the same pair gives, by \(|1\pm e^{it}|^2=2\pm2\cos t\),
\[
 \frac{M_{\ell_\infty^n}(x,y)^2}{\|x\|_\infty^2+\|y\|_\infty^2}
 =\frac12\int_0^{2\pi}(2+2|\cos t|)\,\frac{dt}{2\pi}
 =1+\frac2\pi=B_\infty.
\]
Thus \eqref{eq:finite-p} holds at both endpoints as well.
For finite \(n\geq2\), \(\sin(\pi/n)<\pi/n\), so \(A_n<2\).
Since \(1-2/p>0\) when \(2<p<\infty\), the stated strict
improvement follows.

For \(1<p<2\), put \(\theta=2-2/p\).  Then
\[
 \ell_p^n=[\ell_1^n,\ell_2^n]_\theta,
 \qquad \frac1p=1-\theta+\frac{\theta}{2}.
\]
Using \(\UC(\ell_1^n)=2\) and \(\UC(\ell_2^n)=1\),
Theorem~\ref{thm:interpolation} gives
\[
 \UC(\ell_p^n)\leq2^{1-\theta}1^\theta=2^{2/p-1}.
\]
Conversely, the disjoint coordinate vectors satisfy
\(\|e_1+\zeta e_2\|_p=2^{1/p}\) for \(\zeta\in\T\), and hence
\[
 \UC(\ell_p^n)\geq
 \frac{M_{\ell_p^n}(e_1,e_2)^2}{\|e_1\|_p^2+\|e_2\|_p^2}
 =\frac{2^{2/p}}2=2^{2/p-1}.
\]
The two bounds prove \eqref{eq:finite-low-exact} for \(1<p<2\).
The endpoint values \(2\) and \(1\) at \(p=1\) and \(p=2\), respectively,
were established at the beginning of the proof.
\end{proof}

\begin{remark}
The estimate \eqref{eq:finite-p} separates two effects that the
classical constant does not distinguish: \(B_p\) comes from a
two-coordinate Hadamard configuration, whereas \(A_n\) records how
many distinct phases are available in \(n\) coordinates.  The gap
between them is therefore a concrete finite-dimensional problem.
\end{remark}

\begingroup
\section{Numerical checks}
The following calculations illustrate the preceding results; they do
not replace their proofs. In each example, \(Q_m(x,y)\) denotes the
finite-phase quotient in \eqref{eq:Um} before taking the supremum.

For the Euclidean norm on \(\C^2\), the sum of the \(m\)-th roots of
unity vanishes for \(m\geq2\), so \(Q_m(x,y)=1\) for all nonzero
pairs. For \(x=e_1\), \(y=e_2\), the quotient is \(2\) in
\(\ell_1^2\) and \(1/2\) in \(\ell_\infty^2\), for every \(m\geq2\).

For \(x=(1,1)\) and \(y=(1,-1)\) in \(\ell_\infty^2\), direct
calculation gives
\[
 Q_m(x,y)=1+\frac1m\sum_{k=0}^{m-1}
 \left|\cos\frac{2\pi k}{m}\right|,
 \qquad Q_\infty(x,y)=1+\frac2\pi.
\]
The values and absolute errors, rounded to nine decimal places, are
\begin{center}
\begin{tabular}{@{}rrr@{}}
\toprule
\(m\) & \(Q_m(x,y)\) & \(|Q_m(x,y)-Q_\infty(x,y)|\)\\
\midrule
2  & 2.000000000 & 0.363380228\\
4  & 1.500000000 & 0.136619772\\
8  & 1.603553391 & 0.033066382\\
16 & 1.628417437 & 0.008202336\\
32 & 1.634573149 & 0.002046623\\
64 & 1.636108363 & 0.000511409\\
\bottomrule
\end{tabular}
\end{center}
These are quotients for one fixed pair, not the optimized constants
\(\Um(X)\); their nonmonotonicity does not contradict
Theorem~\ref{thm:refinement}.

For the Hadamard quotient \(B_p\) in \eqref{eq:Bp}, composite
Simpson quadrature on \([0,\pi/2]\) with 200,000 subintervals gives
\begin{center}
\begin{tabular}{@{}rrrrrrr@{}}
\toprule
\(p\) & 1 & 2 & 3 & 4 & 6 & 10\\
\midrule
\(B_p\) & 0.818310 & 1.000000 & 1.127658 & 1.216007 & 1.326755 & 1.435090\\
\bottomrule
\end{tabular}
\end{center}
The values are rounded to six decimal places. They illustrate the
bounds in Theorem~\ref{thm:Lp-high} and do not establish the unknown
extremal value on \(\ell_p^2\) for \(2<p<\infty\).
\endgroup

\section{Open problems for this work}

The following problems remain open.

\begin{enumerate}[label=\textup{(\arabic*)},leftmargin=2.5em]
\item Is the Hadamard lower bound \(B_p\) in \eqref{eq:Bp} the exact
value of \(\UC(\ell_p^2)\) for \(2<p<\infty\)?
\item Characterize the spaces for which
\(\UC(X)=\CNJ(X)\).  Equality holds when \(X\) contains an
isometric copy of \(\ell_1^2\), but fails for every finite
\(\ell_\infty^n\) \rev{with \(n\geq2\)}.
\item Develop \(p\)-circular constants by replacing the square in
\eqref{eq:M-def} with a general exponent, and compare them with
power-type complex convexity and smoothness.
\item Compute the finite constants
\(C_{\mathrm{NJ}}^{[m]}(\ell_\infty^n)\) in the region \(n<m\).
Theorem \ref{thm:limit} provides a certified approximation but not a
closed formula.
\end{enumerate}

\section*{Conclusion}

A single freely chosen complex scalar cannot refine the classical
von Neumann--Jordan constant, because its phase can be absorbed into
the vector variable.  Averaging all phases before optimization avoids
this collapse.  The resulting circular mean has a transparent
Fourier-operator interpretation, a finite-phase bridge back to
Clarkson's constant, \rev{a one-sided isometric Hilbert-space
characterization} and computable values on classical and operator spaces.  The same operator
identifies the invariant with two-term Steinhaus type and makes complex
interpolation available, producing a genuinely dimension-sensitive
\(\ell_p^n\) estimate.  \rev{Comparison with the classical constant gives exact lower values
on the upper half of the \(L_p\) and Schatten scales.}
\rev{In particular, the upper condition \(\UC(X)=1\) already
determines the Hilbert norm, without a separate lower or dual
assumption.}  The quantitative James
inequality turns the sharp circular threshold into an explicit
uniform-nonsquareness gap and hence a fixed point criterion.

The complex convexity comparison is deliberately one-sided: the disc
maximum and the circular \(L_2\)-average are different operations.
Keeping them separate is essential for the validity of the theory.
The remaining two-dimensional \(\ell_p\) extremal problem is stated
explicitly as an open question; no unproved extremal formula is used.

\section*{Acknowledgment}

\noindent H. Li was supported by the National Natural Science Foundation of China (No.12471130).

\end{document}